\documentclass[a4paper,12pt]{article}
\usepackage[T1]{fontenc}
\usepackage[utf8]{inputenc}
\usepackage[english]{babel}
\usepackage{times}
\usepackage{geometry}
\usepackage{afterpage}
\usepackage{booktabs}
\usepackage{tabularx}
\usepackage{epstopdf}
\usepackage{graphicx}
\usepackage{sidecap}
\usepackage{wrapfig}
\usepackage{subfigure}
\usepackage{subcaption}
\usepackage{setspace}
\usepackage{multicol}
\usepackage{multirow}
\usepackage{verbatim}
\usepackage[font=small,labelfont=bf]{caption}
\usepackage{xcolor}
\usepackage{textcomp}
\usepackage[hidelinks]{hyperref}
\usepackage[superscript]{cite}
\usepackage{algorithm}
\usepackage{algorithmicx}
\usepackage{algpseudocode}
\usepackage{float}
\usepackage{dirtytalk}
\usepackage{cancel}
\usepackage{chngpage}
\algdef{SE}[DOWHILE]{Do}{doWhile}{\algorithmicdo}[1]{\algorithmicwhile\ #1}

\usepackage{amsmath}
\usepackage{amsthm}
\usepackage{amssymb}
\usepackage{mathtools}
\newcommand{\numberset}{\mathbb}

\newcommand{\R}{\numberset{R}}

\newcommand{\E}{\numberset{E}}
\newcommand{\g}{\mathbf}

\renewcommand{\epsilon}{\varepsilon}
\renewcommand{\phi}{\varphi}
\DeclarePairedDelimiter{\abs}{\lvert}{\rvert}
\DeclarePairedDelimiter{\norm}{\lVert}{\rVert}

\begin{document}

\begin{Large}
\noindent \textbf{Two new calibration techniques of lumped-parameter mathematical models for the cardiovascular system}
\end{Large}

\noindent \textbf{Andrea Tonini$^1$*,  Francesco Regazzoni$^1$, Matteo Salvador$^2$,  Luca Dede'$^1$, Roberto Scrofani$^3$, Laura Fusini$^4$, Chiara Cogliati$^{5,6}$, Gianluca Pontone$^4$, Christian Vergara$^7$ \& Alfio Quarteroni$^{1,8}$}

\noindent $^1$ MOX, Dipartimento di Matematica, Politecnico di Milano, Milan, Italy.

\noindent $^2$ Institute for Computational and Mathematical Engineering, Stanford University, California, USA.

\noindent $^3$ UOC Cardiochirurgia Fondazione IRCCS Ca' Granda, Ospedale Maggiore Policlinico di Milano, Milan, Italy.

\noindent $^4$ Centro Cardiologico Monzino IRCCS, Milan, Italy.

\noindent $^5$ Internal Medicine, L. Sacco Hospital, Milan, Italy.

\noindent $^6$ Department of Biomedical and Clinical Sciences, Università di Milano, Milan, Italy.

\noindent $^7$ LABS, Dipartimento di Chimica, Materiali e Ingegneria Chimica, Politecnico di Milano, Milan, Italy.

\noindent $^8$ (Professor Emeritus) Institute of Mathematics, Ecole Polytechnique Fédérale de Lausanne, Switzerland

\noindent *email: {\color{blue} \href{mailto:andrea.tonini@polimi.it}{andrea.tonini@polimi.it}}

\section*{Data availability statement}
The data that support the findings of this study are available from the corresponding author upon reasonable request.

\section*{Funding statement}
This work has been funded by the Italian research project FISR (Fondo Integrativo Speciale per la Ricerca) 2020 \say{Mathematical modelling of Covid-19 effects on the cardiac function, Mathematical modelling and analysis of clinical data related to the Covid-19 pandemic in Italy, 2021\_ASSEGNI\_DMAT\_10}. Funding agency: MIUR (Italian Ministry of Education, Universities and Research).

\section*{Conflict of interest disclosure}
The authors declare no competing interests.

\section*{Patient consent statement}
Each patient provided consent to use his/her data for observational studies.

\section*{Acknowledgments}
AT, FR, LD, CV, AQ are members of the INdAM group GNCS \say{Gruppo Nazionale per il Calcolo Scientifico} (National Group for Scientific Computing).

FR has been supported by the INdAM GNCS Project E53C22001930001 and by the project PRIN2022, MUR, Italy, 2023--2025, P2022N5ZNP \say{SIDDMs: shape-informed data-driven models for parametrized PDEs, with application to computational cardiology'}.

CV and LD would like to acknowledge the Italian Ministry of University and Research (MIUR) within the PRIN (Research projects of relevant national interest) MIUR PRIN22-PNRR n. P20223KSS2 \say{Machine learning for fluid-structure interaction in cardiovascular problems: efficient solutions, model reduction, inverse problems}.

CV would like to acknowledge the Italian Ministry of Health within the PNC PROGETTO HUB - DIAGNOSTICA AVANZATA (HLS-DA) \say{INNOVA}, PNC-E3-2022-23683266.

The present research is part of the activities of “Dipartimento di Eccellenza 2023–2027”, MUR, Italy, Dipartimento di Matematica, Politecnico di Milano.

\newpage

\noindent \textbf{Abstract:} Cardiocirculatory mathematical models serve as valuable tools for investigating physiological and pathological conditions of the circulatory system. To investigate the clinical condition of an individual, cardiocirculatory models need to be personalized by means of calibration methods. In this study we propose a new calibration method for a lumped-parameter cardiocirculatory model. This calibration method utilizes the correlation matrix between parameters and model outputs to calibrate the latter according to data. We test this calibration method and its combination with L-BFGS-B (Limited memory Broyden – Fletcher – Goldfarb – Shanno with Bound constraints) comparing them with the performances of L-BFGS-B alone. We show that the correlation matrix calibration method and the combined one effectively reduce the loss function of the associated optimization problem. In the case of in silico generated data, we show that the two new calibration methods are robust with respect to the initial guess of parameters and to the presence of noise in the data. Notably, the correlation matrix calibration method achieves the best results in estimating the parameters in the case of noisy data and it is faster than the combined calibration method and L-BFGS-B. Finally, we present real test case where the two new calibration methods yield results comparable to those obtained using L-BFGS-B in terms of minimizing the loss function and estimating the clinical data. This highlights the effectiveness of the new calibration methods for clinical applications.

\noindent \textbf{Keywords:} Cardiocirculatory models, optimization, parameter estimation, global sensitivity analysis.

\noindent \textbf{Abbreviations:}  \textbf{ODE}: ordinary differential equations; \textbf{CMC}: correlation matrix calibration method; \textbf{L-BFGS-B}: Limited memory Broyden – Fletcher – Goldfarb – Shanno with Bound constraints; \textbf{CMC-L-BFGS-B}: hybrid calibration method between the correlation matrix calibration method and L-BFGS-B; $\mathbf{HR}$: heart rate; $\mathbf{BSA}$: body surface area; $\mathbf{LA_{\mathrm{Vmax}}}$: maximal left atrial volume; $\mathbf{LV_{\mathrm{EDV}}}$: left ventricular end diastolic volume; $\mathbf{LV_{\mathrm{ESV}}}$: left ventricular end systolic volume; $\mathbf{LV_{\mathrm{EF}}}$: left ventricular ejection fraction; $\mathbf{max \nabla P_{\mathrm{rAV}}}$: maximal right atrioventricular pressure gradient; $\mathbf{SAP_{\mathrm{max}}}$: systolic systemic arterial pressure; $\mathbf{SAP_{\mathrm{min}}}$: diastolic systemic arterial pressure; $\mathbf{PAP_{\mathrm{max}}}$: systolic pulmonary arterial pressure; $\mathbf{LA_\mathrm {Pmax}}$: maximal left atrial pressure; $\mathbf{LA_\mathrm {Pmin}}$: minimal left atrial pressure; $\mathbf{LA_\mathrm {Pmean}}$: mean left atrial pressure; $\mathbf{LV_{\mathrm{SV}}}$: left ventricular stroke volume; $\mathbf{CO}$: cardiac output; $\mathbf{CI}$: cardiac index;   $\mathbf{LV_{\mathrm{Pmax}}}$: maximal left ventricular pressure; $\mathbf{LV_{\mathrm{Pmin}}}$: minimal left ventricular pressure; $\mathbf{RA_{\mathrm{Vmax}}}$: maximal right atrial volume; $\mathbf{RA_\mathrm{Pmax}}$: maximal right atrial pressure; $\mathbf{RA_\mathrm{Pmin}}$: minimal right atrial pressure; $\mathbf{RA_\mathrm{Pmean}}$: mean right atrial pressure; $\mathbf{RV_{\mathrm{EDV}}}$: right ventricular end diastolic volume; $\mathbf{RV_{\mathrm{ESV}}}$: right ventricular end systolic volume; $\mathbf{RV_{\mathrm{EF}}}$: right ventricular ejection fraction;  $\mathbf{RV_{\mathrm{Pmax}}}$: maximal right ventricular pressure; $\mathbf{RV_{\mathrm{Pmin}}}$: minimal right ventricular pressure; $\mathbf{PAP_{\mathrm{min}}}$: diastolic pulmonary arterial pressure; $\mathbf{PAP_{\mathrm{mean}}}$: mean pulmonary arterial pressure; $\mathbf{PWP_{min}}$: minimal pulmonary wedge pressure; $\mathbf{PWP_{mean}}$: mean pulmonary wedge pressure; $\mathbf{SVR}$: systemic vascular resistance; $\mathbf{PVR}$: pulmonary vascular resistance 

\section{Introduction}

Cardiocirculatory mathematical models have been developed to reproduce physiological and pathological conditions of the human body\cite{shi2011review, quarteroni2016geometric,dede2021modeling,fernandes2021integrated,de2006modelling,shi2006numerical, tonini2024mathematical,zhu2023computational}. \emph{Lumped-parameter} models (named also 0D models) are set on a partition of the cardiovascular system into different reduced compartments (e.g., systemic arteries or left atrium), where only the average flow rates and pressures are computed at each time. Each compartment is characterized by a set of parameters (e.g., resistances of the vessels or elastances of the cardiac chambers) that typically refer to average, physiological conditions. 0D models 

Given a set of either in silico generated or clinical data, a calibration method modifies the parameters of the lumped-parameter model to minimize the distance between data and model outputs thus making the model patient-specific \cite{dede2021modeling,herrmann2020modeling}. Previous works studied the advantages of calibration methods for cardiocirculatory lumped-parameter models where synthetically generated data are used in the minimization process. For example, Laubscher et al.\cite{laubscher2023estimation} used a combination of Adam optimizer and L-BFGS-B (Limited memory Broyden – Fletcher – Goldfarb – Shanno with Bound constraints) to estimate different parameters that significantly influence the left ventricular pressure-volume loop; Bjørdalsbakke et al.\cite{bjordalsbakke2022parameter} employed a trust region reflective algorithm in ten different estimation scenarios following a sensitivity analysis. These works used gradient-based calibration methods. Saxton et al.\cite{saxton2023personalised} used an Unscented Kalman filter to estimate the parameters of a cardiocirculatory model that accounts for the left ventricle and for the systemic circulation with a time varying heartbeat period.

In this work, we propose a new gradient-free method with the aim of providing an efficient calibration procedure in the context of lumped-parameter cardiocirculatory models. This method uses of the correlation matrix between parameters and model outputs to surrogate the gradient of the loss function. We refer to this method as \emph{Correlation Matrix Calibration} (CMC) method. We identify the parameters that significantly affect a set of model outputs corresponding to a set of available data (such as the maximal left atrial volume or the systemic systolic arterial pressure) by means of a global sensitivity analysis \cite{saltelli2002making,sobol1993sensitivity} and we perform the calibration procedure of the identified parameters.

We compare the proposed CMC method to the L-BFGS-B method\cite{byrd1995limited,zhu1997algorithm} on a dataset of in silico generated data. Moreover, we also propose a combination of the two calibration methods (CMC-L-BFGS-B) in order to achieve better results in terms of accuracy with respect to CMC alone. Specifically, we first apply CMC and then L-BFGS-B, because the first one acts as a global method, whereas the latter as a local method that can improve the accuracy of the former. Although the proposed new methods could be applied to models different from the cardiocirculatory ones, we restric ourselves only to the latter case.

To further validate CMC and CMC-L-BFGS-B, we apply them to patient-specific data related to COVID-19 pneumonia provided by Centro Cardiologico Monzino and L. Sacco Hospital in Milan, Italy. In severe COVID-19-related pneumonia, right ventricle involvement seems to mainly drive cardiac function damages, while consequences on the left ventricle appear to be less common \cite{dandel2022heart}. Right ventricle dilation, diminished right ventricular function and elevated pulmonary arterial systolic pressure are associated with mortality in severe COVID-19 \cite{park2020eye,diaz2021association}. Moreover, endothelial damages with diffuse micro-thrombosis has been widely described in histological studies in COVID-19 pneumonia patients causing an increase in pulmonary resistances and a reduction in pulmonary compliances \cite{mauri2020potential,peng2020using}. We choose the ranges where the parameters vary during the calibration procedure according to these observations.

The outline of this paper is as follows. In Section \ref{sec:model}, we describe the lumped-parameter cardiocirculatory model together with its parameters and the computable model outputs. In Section \ref{sec:newcal}, we present the global sensitivity analysis and the considered calibration methods. In Section \ref{sec:results}, we test the robustness of the calibration methods both on in silico generated data and on patient-specific data. In Section \ref{sec:concl}, we draw the conclusions of this work.

\section{Lumped-parameter cardiocirculatory model}\label{sec:model}
In this section we introduce the lumped-parameter cardiocirculatory model togheter with its parameters and its outputs.

A lumped-parameter cardiocirculatory model describes the human cardiovascular system as an electrical circuit: the current represents the blood flow through vessels and valves, the electric potential corresponds to the blood pressure, the electric resistance plays the role of the resistance to blood flow, the capacitance represents the vessel compliance and the inductance corresponds to the blood inertia.

A lumped-parameter cardiocirculatory model partitions the cardiovascular system into distinct compartments (e.g. right atrium, systemic arteries/veins). For each of them, a system of ordinary differential equations (ODE) \cite{shi2011review,quarteroni2016geometric,dede2021modeling} describes the time evolution of the unknowns of the model (pressures, flow rates and cardiac volumes).
 
\begin{figure}[H]
\centering
\includegraphics[width=0.5\linewidth]{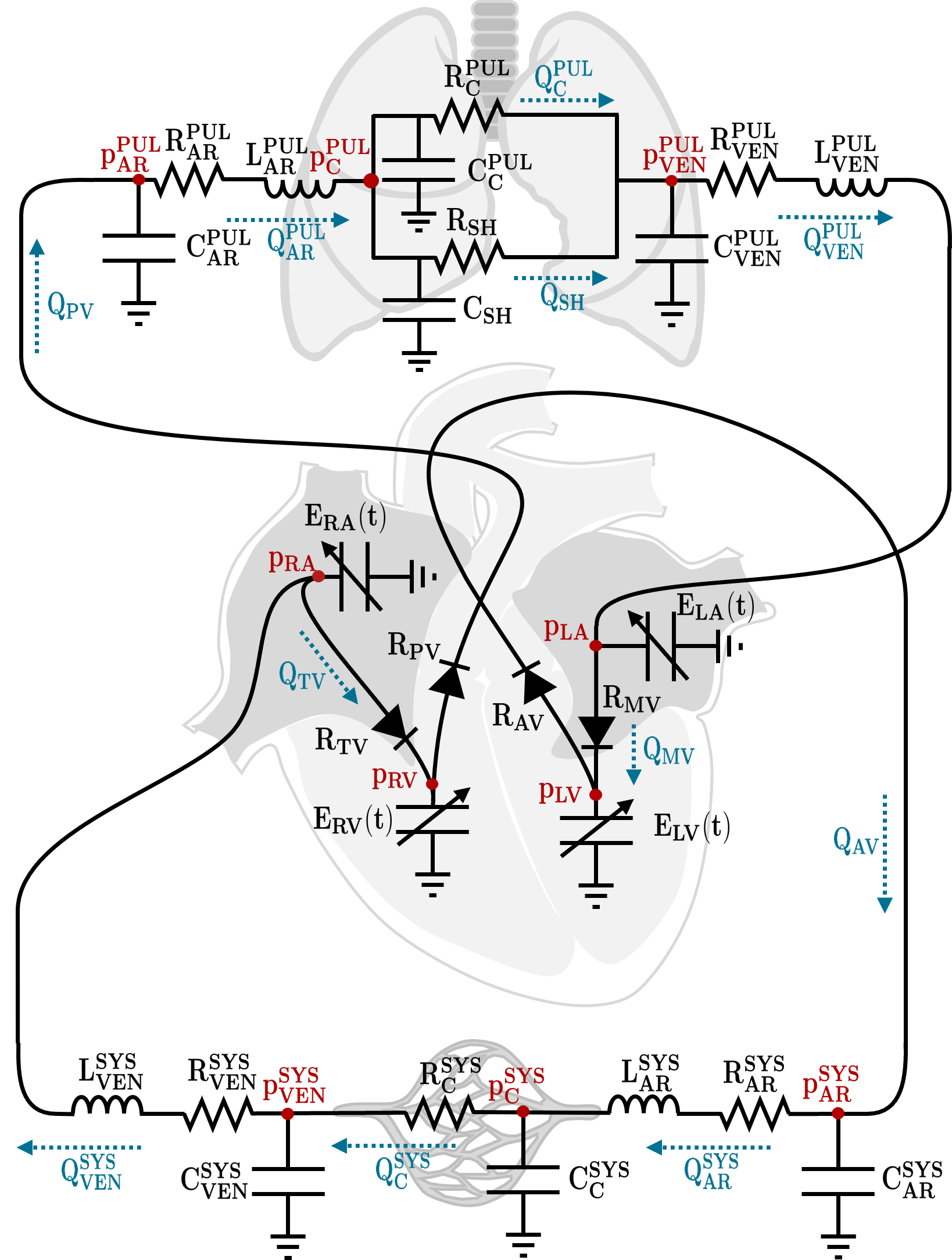}
\caption{Lumped-parameter model. We depict pressures and flow rates in red and blue, respectively, and parameters in black.}
\label{fig:lumpedparametercardiovascular}
\end{figure}

The proposed lumped-parameter model is a modification of the one proposed by Regazzoni et al. \cite{regazzoni2022cardiac}, that we improved by adding further compartments representing the systemic and pulmonary micro-vasculature (Figure \ref{fig:lumpedparametercardiovascular}). The model consists in the four cardiac chambers, the systemic and pulmonary circulation, split into arterial, capillary and venous compartments. We describe the pulmonary microcirculation by making use of two compartments accounting for the oxygenated and non-oxygenated capillaries, respectively, because some pulmonary capillaries do not oxygenate due to hypoxic vasoconstriction in regions with a low oxygen concentration in the alveoli, even if blood perfusion of lungs is normal. The fraction of blood that does not oxygenate is called pulmonary shunt. In healthy conditions pulmonary shunt is lower than $5\%$ \cite{velthuis2015pulmonary}, whereas in presence of acute respiratory distress sindromes, as COVID-19, it can increases up to $60\%$ \cite{saha2021correlation,gattinoni2020covid}. The hemodynamic of the cardiovascular system is described by means of a dynamical system:
\begin{align}
\begin{cases}
\dot {\mathbf x}(t;\mathbf x_\mathrm{0},\mathbf p) = \mathbf f(t,\mathbf x(t;\mathbf x_\mathrm{0},\mathbf p);\mathbf p)\qquad t \in (0,T] \\
\mathbf x(0;\mathbf x_\mathrm{0},\mathbf p) = \mathbf x_\mathrm{0}
\end{cases}
\end{align}
where $\mathbf x$, $\mathbf x_\mathrm{0}$, $\mathbf p$, $\mathbf f$ and $T$ represent the state variables, the initial conditions, the parameters, the system right hand side (rhs) and the final time, respectively. Moreover, we can calculate some additional quantities $\mathbf y$, that we call model outputs, as functions (e.g. the maximum or the mean) of the state variables $\mathbf x$ and the parameters $\mathbf p$:
\[
\mathbf y = \mathbf g(\mathbf x(t;\mathbf x_\mathrm{0}, \mathbf p);\mathbf p).
\]

The state variables $\mathbf x$ of our model are the volumes of the left atrium ($V_\mathrm{LA}$) and ventricle ($V_\mathrm{LV}$) and of the right atrium ($V_\mathrm{RA}$) and ventricle ($V_\mathrm{RV}$), the circulatory pressures of the systemic arteries ($p_\mathrm{AR}^\mathrm{SYS}$), capillaries ($p_\mathrm{C}^\mathrm{SYS}$) and veins ($p_\mathrm{VEN}^\mathrm{SYS}$) and of the pulmonary arteries ($p_\mathrm{AR}^\mathrm{PUL}$), capillaries ($p_\mathrm{C}^\mathrm{PUL}$) and veins ($p_\mathrm{VEN}^\mathrm{PUL}$), the circulatory fluxes of the systemic arteries ($Q_\mathrm{AR}^\mathrm{SYS}$) and veins ($Q_\mathrm{VEN}^\mathrm{SYS}$) and of the pulmonary arteries ($Q_\mathrm{AR}^\mathrm{PUL}$) and veins ($Q_\mathrm{VEN}^\mathrm{PUL}$).

The lumped-parameter cardiocirculatory model depends on the heart rate (HR), that determines the heartbeat period T$_{\text{HB}}$$=60/$HR, and on the parameters $\mathbf p$ reported in Table \ref{table:params}. We fix the final time $T = 25 T_\mathrm{\text{HB}}$ to reach the limit cycle of the dynamical system.

The dynamics of cardiac blood volumes accounts for the inward and outward fluxes:
\begin{gather}
\dot V_{\mathrm{RA}}(t) = Q_{\mathrm{VEN}}^{\mathrm{SYS}}(t)-Q_{\mathrm{TV}}(t), \qquad \dot V_{\mathrm{LA}}(t) = Q_{\mathrm{VEN}}^{\mathrm{PUL}}(t)-Q_{\mathrm{MV}}(t),\\
\dot V_{\mathrm{RV}}(t) = Q_{\mathrm{TV}}(t)-Q_{\mathrm{PV}}(t), \qquad \dot V_{\mathrm{LV}}(t) = Q_{\mathrm{MV}}(t)-Q_{\mathrm{AV}}(t),
\end{gather}
where $Q_\mathrm{MV},Q_\mathrm{AV},Q_\mathrm{TV}$ and $Q_\mathrm{PV}$ are the blood flows through the mitral, aortic, tricuspid and pulmonary valves, respectively. These flows depend on the pressure jump from the upstream to the downstream compartment:
\begin{gather}
Q_{\mathrm{TV}}(t) = Q_{\mathrm{valve}}(p_{\mathrm{RA}}(t)-p_{\mathrm{RV}}(t)), \quad Q_{\mathrm{MV}}(t) = Q_{\mathrm{valve}}(p_{\mathrm{LA}}(t)-p_{\mathrm{LV}}(t)),\\
Q_{\mathrm{PV}}(t) = Q_{\mathrm{valve}}(p_{\mathrm{RV}}(t)-p_{\mathrm{AR}}^{\mathrm{PUL}}(t)), \quad Q_{\mathrm{AV}}(t) = Q_{\mathrm{valve}}(p_{\mathrm{LV}}(t)-p_{\mathrm{AR}}^{\mathrm{SYS}}(t)),
\end{gather}
where $p_\mathrm{LA}$, $p_\mathrm{LV}$, $p_\mathrm{RA}$ and $p_\mathrm{RV}$ are the unknown pressures inside the left atrium, the left ventricle, the right atrium and the right ventricle, respectively. $Q_\mathrm{valve}$ is the blood flow across a valve that depends on the pressure jump across the valve:
\begin{align}
Q_{\mathrm{valve}}(\Delta p) = \frac{\Delta p}{R_{\mathrm{valve}}(\Delta p)}.
\end{align}

The resistance of the leaflets of each valve is
\begin{align}
R_{\mathrm{valve}}(\Delta p) = \sqrt{R_{\mathrm{min}}R_{\mathrm{max}}} \left( \frac{R_{\mathrm{max}}}{R_{\mathrm{min}}} \right)^\mathrm{\frac{\text{atan} \left(-100\pi \Delta p\right)}{\pi}}.
\end{align}
$R_{\mathrm{valve}}$ ranges from $R_{\mathrm{min}}$ (with $\Delta p \to +\infty$) to $R_{\mathrm{max}}$ (with $\Delta p \to -\infty$), where $R_{\mathrm{min}}$ and $R_{\mathrm{max}}$ are the minimal (open valve) and maximal (closed valve) resistances given by the leaflets of the valves.

We model each cardiac chamber as a pressure generator. The pressure generated  $p_{\mathrm{c}}(t)$ by the myocardium of the cardiac chamber $c$ ($c \in \{LA,LV,RA,RV\}$) depends on the blood volume contained in the chamber itself $V_{\mathrm{c}}(t)$, the unloaded volume $V_{U,\mathrm{c}}$ (i.e. volume at zero pressure) and the time varying elastance $E_{\mathrm{c}}(t)$, that represents the cardiac chamber contractility. 
\begin{gather}
p_{\mathrm{c}}(t) = E_{\mathrm{c}}(t)(V_{\mathrm{c}}(t)-V_{U,\mathrm{c}}).
\end{gather}

The time varying elastance depends on the passive elastance $EB_{\mathrm{c}}$ (i.e. the inverse of the cardiac chamber compliance), the maximum active elastance $EA_{\mathrm{c}}$ and a periodic function $e_{\mathrm{c}}(t)$ that accounts for cardiac activation phases\cite{dede2021modeling,liang2009multi}.
\begin{gather}
E_{\mathrm{c}}(t) = EB_{\mathrm{c}}+EA_{\mathrm{c}}e_{\mathrm{c}}(t),\\
{\scriptsize
e_{\mathrm{c}}(t) =
\begin{cases}
\frac{1}{2}\left[1-\cos\left(\frac{\pi}{TC_{\mathrm{c}}}\mathrm{mod}\left(t-tC_{\mathrm{c}},T_{\mathrm{HB}} \right) \right)\right] \qquad &\text{if } 0 \le \mathrm{mod}\left(t-tC_{\mathrm{c}},T_{\mathrm{HB}}\right)<TC_{\mathrm{c}},\\[5pt]
\frac{1}{2}\left[1+\cos\left(\frac{\pi}{TR_{\mathrm{c}}}\mathrm{mod}\left(t-tR_{\mathrm{c}},T_{\mathrm{HB}} \right) \right)\right] \qquad &\text{if } 0 \le \mathrm{mod}\left(t-tR_{\mathrm{c}},T_{\mathrm{HB}}\right)<TR_{\mathrm{c}},\\[5pt]
0 & \mathrm{otherwise},
\end{cases}\label{eq:el}
}
\end{gather}
where $TC_{\mathrm{c}}$, $TR_{\mathrm{c}}$, $tC_{\mathrm{c}}$ and $tR_{\mathrm{c}}$ are the durations of the contraction and relaxation phases and the times of the contraction and relaxation phases of the cardiac chamber, respectively (Table \ref{table:params}). The first and the second equations of (\ref{eq:el}) account for the contraction and relaxation phases of the cardiac chamber, respectively.

\begin{figure}[H]
	\centering
 	\subfigure[]{\includegraphics[width=0.47\linewidth]{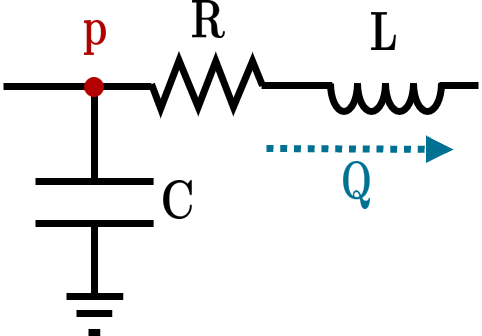}} 
	\subfigure[]{\includegraphics[width=0.47\linewidth]{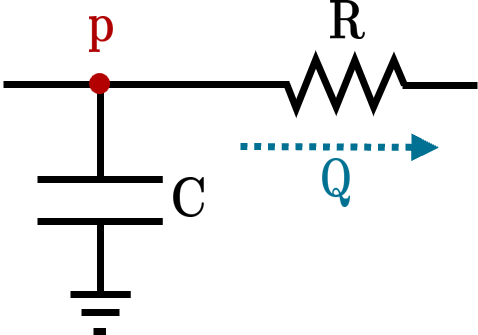}}
	\caption{RLC Windkessel circuit used for the arterial and venous compartments (a) and RC Windkessel circuit used for capillary compartments (b).}
	\label{fig:windkessel}
\end{figure}

Each circulatory compartment (e.g. systemic arterial circulation) is modeled as a Windkessel circuit (Figure \ref{fig:windkessel}) describing the circulation by means of the Kirchhoff's circuit laws. For example, for the systemic arterial circulation we have: 
\begin{gather}
C_{\mathrm{AR}}^{\mathrm{SYS}}\dot p_{\mathrm{AR}}^{\mathrm{SYS}}(t) = Q_{\mathrm{AV}}(t)-Q_{\mathrm{AR}}^{\mathrm{SYS}}(t),\label{eq:pressure}\\
L_{\mathrm{AR}}^{\mathrm{SYS}}\dot Q_{\mathrm{AR}}^{\mathrm{SYS}}(t) = -Q_{\mathrm{AR}}^{\mathrm{SYS}}(t)R_{\mathrm{AR}}^{\mathrm{SYS}}+p_{\mathrm{AR}}^{\mathrm{SYS}}(t)-p_\mathrm{C}^{\mathrm{SYS}}(t), \label{eq:flux}
\end{gather}
where $R_\mathrm{AR}^{SYS}$, $C_\mathrm{AR}^{SYS}$ and $L_\mathrm{AR}^{SYS}$ are the systemic arterial resistance, compliance and inertia, respectively (Table \ref{table:params}). Analogous equations hold for the systemic venous circulation and for the pulmonary arterial and venous circulation. Notice that for the capillary circulation the blood inertia is negligible \cite{albanese2016integrated}, so we set it to $0$. Summarizing, the whole dynamical system is:
\[
\begin{cases}
\dot V_\mathrm{LA}(t)=Q_\mathrm{VEN}^\mathrm{PUL}(t)-Q_\mathrm{MV}(t)\\[2pt] 
\dot V_\mathrm{LV}(t)=Q_\mathrm{MV}(t)-Q_\mathrm{AV}(t)\\[2pt]
 C_\mathrm{AR}^\mathrm{SYS}\dot p_\mathrm{AR}^\mathrm{SYS}(t)=Q_\mathrm{AV}(t)-Q_\mathrm{AR}^\mathrm{SYS}(t)\\[2pt]
L_\mathrm{AR}^\mathrm{SYS}\dot Q_\mathrm{AR}^\mathrm{SYS}(t)=-R_\mathrm{AR}^\mathrm{SYS}Q_\mathrm{AR}^\mathrm{SYS}(t)+p_\mathrm{AR}^\mathrm{SYS}(t)-p_\mathrm{C}^\mathrm{SYS}(t)\\[2pt]
C_\mathrm{C}^\mathrm{SYS}\dot p_\mathrm{C}^\mathrm{SYS}(t) = Q_\mathrm{AR}^\mathrm{SYS}(t)-Q_\mathrm{C}^\mathrm{SYS}(t)\\[2pt]
C_\mathrm{VEN}^\mathrm{SYS}\dot p_\mathrm{VEN}^\mathrm{SYS}(t)=Q_\mathrm{C}^\mathrm{SYS}(t)-Q_\mathrm{VEN}^\mathrm{SYS}(t)\\[2pt]
L_\mathrm{VEN}^\mathrm{SYS} \dot Q_\mathrm{VEN}^\mathrm{SYS}(t)=-R_\mathrm{VEN}^\mathrm{SYS}Q_\mathrm{VEN}^\mathrm{SYS}(t)+p_\mathrm{VEN}^\mathrm{SYS}(t)-p_\mathrm{RA}(t)\\[2pt]
\dot V_\mathrm{RA}(t)=Q_\mathrm{VEN}^\mathrm{SYS}(t)-Q_\mathrm{TV}(t)\\[2pt]
\dot V_\mathrm{RV}(t)=Q_\mathrm{TV}(t)-Q_\mathrm{PV}(t)\\[2pt]
C_\mathrm{AR}^\mathrm{PUL}\dot p_\mathrm{AR}^\mathrm{PUL}(t)=Q_\mathrm{PV}(t)-Q_\mathrm{AR}^\mathrm{PUL}(t)\\[2pt]
L_\mathrm{AR}^\mathrm{PUL}\dot Q_\mathrm{AR}^\mathrm{PUL}(t)=-R_\mathrm{AR}^\mathrm{PUL}Q_\mathrm{AR}^\mathrm{PUL}(t)+p_\mathrm{AR}^\mathrm{PUL}(t)-p_\mathrm{C}^\mathrm{PUL}(t)\\[2pt]
(C_\mathrm{SH}+C_\mathrm{C}^\mathrm{PUL})\dot p_\mathrm{C}^\mathrm{PUL}(t) = Q_\mathrm{AR}^\mathrm{PUL}(t)-Q_\mathrm{SH}(t)-Q_\mathrm{C}^\mathrm{PUL}(t)\\[2pt]
C_\mathrm{VEN}^\mathrm{PUL}\dot p_\mathrm{VEN}^\mathrm{PUL}(t)=Q_\mathrm{SH}(t)+Q_\mathrm{C}^\mathrm{PUL}(t)-Q_\mathrm{VEN}^\mathrm{PUL}(t)\\[2pt]
L_\mathrm{VEN}^\mathrm{PUL}\dot Q_\mathrm{VEN}^\mathrm{PUL}(t)=-R_\mathrm{VEN}^\mathrm{PUL}Q_\mathrm{VEN}^\mathrm{PUL}(t)+p_\mathrm{VEN}^\mathrm{PUL}(t)-p_\mathrm{LA}(t)\\[2pt]
\end{cases}
\]
coupled with suitable initial conditions.

Once we solve the discretized dynamical system, finding an approximation of the state variables, we compute the model outputs $\mathbf y$. We define the model outputs in Table \ref{table:qoi}, where a subscript “I-” refers to the indexed volumes, i.e. normalized by the body surface area (BSA). The calibration procedure will be performed over specific model outputs, in particular those reported in the top part of Table \ref{table:qoi}.

We determine the reference setting of parameters $\mathbf p^\mathrm{R}$ to reproduce an ideal healthy individual (see Table \ref{table:params}). Specifically, we fix HR$=80$ bpm \cite{zhang2009heart} and we find the values of the other parameters as a modification of literature values \cite{albanese2016integrated,dede2021modeling}, in such a way that the model outputs $\mathbf y$ lie in echocardiographic ranges related to a healthy individual (Table \ref{table:qoi}).

\section{Novel calibration methods} \label{sec:newcal}

In this section, we first describe the global sensitivity analysis toghether with the general aim of the calibrations (Section \ref{sec:GSA}). This procedure is common to all the calibration methods presented in what follows and allows to restrict the calibration procedure to a lower number of parameters. Then, we describe the new calibration method in Section \ref{sec:cor} and an hybrid method that combine CMC and L-BFGS-B in Section \ref{sec:CMC-LBFGS_B} (CMC-L-BFGS-B).

\subsection{Global sensitivity analysis and aim of the calibration} \label{sec:GSA}
We perform a global sensitivity analysis to determine which parameters affect significantly the model outputs related to the data (Table \ref{table:qoi}). Letting the parameters vary randomly in a hyperbox, we estimate total Sobol indices \cite{sobol1993sensitivity}, that evaluate the impact of a parameter $p_\mathrm{k}$ on a certain model output $y_\mathrm{j}$, accounting also for high-order interactions among parameters:
\[
S_\mathrm{k}^\mathrm{j,T}=1-\frac{Var_\mathrm{\mathbf p_\mathrm{\sim k}}[\numberset{E}_\mathrm{p_\mathrm{k}}[y_\mathrm{j}|\mathbf p_\mathrm{\sim k}]]}{Var[y_\mathrm{j}]}
\]
where $\mathbf p_\mathrm{\sim k}$ indicates the set of all parameters excluding the $k^\mathrm{th}$ one. $\numberset{E}$ and $Var$ are the expected value and the variance, respectively, and the subscripts indicate the random variable measure to use for the integration. When the subscript is absent, the integration is performed with respect to all the random variables. $\numberset{E}_\mathrm{p_\mathrm{k}}[y_\mathrm{j}|\mathbf p_\mathrm{\sim k}]$ is the expected value of the model output $y_\mathrm{j}$ conditioned with respect to the parameters $\mathbf p_\mathrm{\sim k}$.

We estimate total Sobol indices by sampling the parameters in a hyperbox employing the Saltelli's method \cite{saltelli2002making} and we compute the corresponding model outputs by means of the lumped-parameter cardiocirculatory model. We build the hyperbox around the ideal healthy reference setting of parameters $\mathbf{p}^\mathrm{R}$ (Appendix \ref{app:bounds}). We do not compute the sensitivity to the $HR$, as this parameter is easily measured from patients, so we are not interested in its calibration from indirect measures. The Saltelli's method allows for a linear increase in the number of samples with respect to the number of varying parameters $N_\mathrm{p}=32$: the number of samples is $2N(N_\mathrm{p}+1)$ where $N$ is a user defined variable. We choose $N=1024$ that corresponds to $67584$ samples and allows for small confidence intervals of the total Sobol indices.

We stress that the shape and the position of the hyperbox where we sample the parameters affect the Sobol indices. In particular, the wider the interval for a parameter is, the higher the Sobol indices for that parameter will be, because it can potentially generate higher variability in the model outputs. Therefore, the Sobol indices computed here are specific for the chosen ranges of the parameters.

According to the values of the total Sobol indices we select the parameters to calibrate.

During the calibration procedure (based on data generated in silico or measured from specific patients), the parameters eligible for the calibration procedure can vary in the ranges described in Appendix \ref{app:bounds}. Even if we fixed the ranges of the parameters, they are large enough to cover a wide range of possible conditions.

The goal of the calibration methods is to minimize the loss function (mean squared error):
\begin{gather}
\label{eq:loss}
\mathcal MSE(\mathbf p) = \frac{1}{N_\mathrm{d}}\sum_{\mathrm{i=1}}^{\mathrm{N_d}}\left( \frac{d_\mathrm{i}-y_\mathrm{j(i)} (\mathbf p)}{d_\mathrm{i}} \right)^\mathrm{2}
\end{gather}
where $N_\mathrm{d}$ is the number of data $d_\mathrm{i}$, for $1\le i \le N_d$, at disposal. Each data is represented by a model output. We indicate with $j(i)$ the index of the entry of the model outputs vector $\mathbf y$ approximating the $i$-th data. $\mathbf p$ is the set of parameters. 

\subsection{Correlation matrix calibration method}\label{sec:cor}

CMC uses a surrogate of the gradient of the loss function. The Pearson correlation coefficient measures the linear relationship between two variables \cite{pearson1895vii}. If $f:\R \to \R$ is a (non constant) function such that $f(x) \le f(y)$ for almost every $x<y$ with $x,y \in \R$ and $x,y$ are sampled from an absolute continuous distribution $X \sim g$ then the Pearson's correlation coefficient is non-negative: 
\[
\rho_{X,f(X)} = \frac{\E[(X-\E[X])(f(X)-\E[f(X)])]}{\sqrt{\E[(X-\E[X])^\mathrm{2}]\; \E[(f(X)-\E[f(X)])^\mathrm{2}]}}\ge 0
\]
Indeed, the numerator determines the sign of the correlation coefficient:
\begin{align*}
&\E[(X-\E[X])(f(X)-\E[f(X)])] \\
& =\E[(X-\E[X])(f(X)-f(\E[X])+f(\E[X])-\E[f(X)])]\\
& =\E[(X-\E[X])(f(X)-f(\E[X]))] + (f(\E[X])-\E[f(X)])\E[X-\E[X]]\\
& = \E[(X-\E[X])(f(X)-f(\E[X]))]\\
& = \int_{-\infty}^{+\infty}(x-\E[X])(f(x)-f(\E[X]))g(x)dx
\end{align*}
Using the monotonicity of $f$ almost everywhere, the signs of $x-\E[X]$ and $f(x)-f(\E[X])$ are concordant and so $\rho_{X,f(X)}\ge 0$. The converse holds if $f$ is monotone non-increasing. Given $n$ paired data points $\{(x_\mathrm{i},f(x_\mathrm{i}))\}_\mathrm{i=1}^\mathrm{n}$ with $\overline{f(x)} = \frac{1}{n}\sum_\mathrm{i=1}^\mathrm{n}f(x_\mathrm{i}) $, a similar argument applies to the sample Pearson's correlation coefficient\cite{pearson1895vii}
\[
r_{X,f(X)} = \frac{\sum_\mathrm{i=1}^\mathrm{n}(x_\mathrm{i}-\bar x)(f(x_\mathrm{i})-\overline{f(x)})}{\sqrt{ \left (\sum_\mathrm{i=1}^\mathrm{n}(x_\mathrm{i}-\bar x)^\mathrm{2} \right ) \left (\sum_\mathrm{i=1}^\mathrm{n}(f(x_\mathrm{i})-\overline{f(x)})^\mathrm{2} \right )}}
\]
with the assumption of $f$ being monotone everywhere. Therefore, we use the sample correlation coefficient between parameters and model outputs to get information about the monotonicity of the loss function and its gradient with respect to the parameters.

To compute the correlation matrix $M$ (where $M_\mathrm{l,j}=r_\mathrm{p_\mathrm l,y_\mathrm j}$), we perform $3200$ tests ($100N_\mathrm{p}$ where $N_\mathrm p = 32$ in our case). We sample the parameters from a uniform distribution in the previously mentioned hyperbox and we compute the corresponding model outputs. The correlation coefficients are local information of the linear relationships between parameters and model outputs around the mean of the parameters, that coincides with the centre of the hyperbox because we sample the parameters from a uniform distibution. Therefore, CMC could lose in accuracy when the setting of parameters is far from the centre of the hyperbox during the calibration procedure. Nonetheless, it achieves good results for the cardiocirculatory model, as we show in Section \ref{sec:results}.

CMC for the cardiovascular system relies on Algorithm \ref{alg:calibration} described in what follows:
\begin{enumerate}
\item Initialize the variables used troughout the calibration procedure (lines $1$-$3$);
\item Run a model simulation and compute the loss function (lines $5$-$6$); \label{item:simulation}
\item If the loss function is greater than a threshold compute the relative errors between data and the related model outptus. Check if there exists a model output that during the current iteration of the calibration procedure has not yet been considered to calibrate the parameters (lines $7$-$11$);
\item Focus on the model output associated to the largest error (line $12$);
\item Consider the parameter associated to the largest correlation coefficient in absolute value related to this model output (line $13$-$14$); \label{item:condition}
\item Avoid considering always the same parameter for the calibration procedure at each iteration (lines $15$-$18$);
\item If a parameter satisfying the condition at step \ref{item:condition} and associated to a correlation coefficient greater than $0.05$ in absolute value does not exist, go on to the next model output associated to the highest error (lines $19$-$28$) and repeat from step \ref{item:condition}. If not, surrogate the derivative of the loss function with respect to the chosen parameter using the correlation coefficients between the parameter and all the model outputs related to data (lines $29$-$40$);
\item If the surrogate gradient of the loss function is less than $0$, perform a calibration step and repeat the whole procedure from step \ref{item:simulation}. If not, go on to the next data (lines $41$-$50$).
\end{enumerate}

\begin{algorithm}[H]
	\captionsetup{labelfont={sc,bf}, labelsep=newline}
	\caption{CMC algorithm}
	\label{alg:calibration}
	\begin{algorithmic}[1]
	\small
		\State Initialize: $\g d$, $tol$, $it_\mathrm{max}$, $it = 0$, $m = \varnothing$
		\Statex \Comment{Data, tolerance, max number of iterations, current iteration, set of calibrated parameters}
		\State Load $M$ \Comment{Load the correlation matrix}
		\State Choose the initial guess of parameters $\g p$ \Comment{$\g p$ lies in the hyperbox used to build $M$}
		\Do \label{algstep:it_condition}
			\State Run a model simulation
			\State Compute $\mathcal L(\g p)$
			\If{$\mathcal L(\g p)\ge tol$}
				\State $e_\mathrm i =  \frac{d_\mathrm{i}-y_\mathrm{j(i)} (\mathbf p)}{d_\mathrm{i}} \qquad \text{for } i = 1,\dots,N_\mathrm d$
				\State $\g a = \g e$ \Comment{Auxiliary vector}
				\State $par_\mathrm{con} = 0$ \Comment{Number of considered parameters for a specific data}
				\While{$\norm{\g a} > 0$} \label{alg:err_condition}
					\State $\bar \imath = \text{argmax}_\mathrm{i = 1,\dots,N_d}\abs{a_\mathrm i}$
					\State $b = M_\mathrm{:,j(\bar \imath)}$ \Comment{Correlation coefficients between $y_\mathrm{j(\bar \imath)}$ and parameters}
					\State  $\bar l = \text{argmax}_\mathrm{l = 1,\dots,N_p}\left| b_\mathrm{l} \right|$
					\While{$p_\mathrm{\bar l} \in m$ \& $\abs{b_\mathrm{\bar l}> 0.05}$} \Comment{Avoid to calibrate the same parameter}
						\State $b_\mathrm{\bar l} = 0$ 
						\State $\bar l = \text{argmax}_\mathrm{l = 1,\dots,N_p}\left| b_{\mathrm {l}} \right|$
					\EndWhile
\algstore{alg:calibration}
\end{algorithmic}
\end{algorithm}

\begin{algorithm}[H]
\begin{algorithmic}
\small
\algrestore{alg:calibration}
					\If{$\abs{b_\mathrm{\bar l}}\le 0.05$} \Comment{Avoid to consider a parameter not significant for $\mathcal L$}
						\If{$\abs m = par_\mathrm{con}$} \Comment{All parameters were used to calibrate only $y_\mathrm{j(\bar \imath)}$}
							\State $a_\mathrm{\bar \imath} = 0$ \Comment{Change error component to decrease}
						\Else \Comment{Some parameters were used to calibrate $y_\mathrm{j(i)}$ with $i \not = \bar \imath$}
							\State $b = M_\mathrm{:, j(\bar \imath)}$
						\EndIf
						\State $par_\mathrm{con} = 0$
						\State $m = \varnothing$
						\State Go to step \ref{alg:err_condition}
					\EndIf
					\State $m = m \cup \{\bar l\}$
					\State $par_\mathrm{con} = par_\mathrm{con}+1$ 
					\State $\hat \nabla \mathcal L = 0$ \Comment{Surrogate $(\nabla \mathcal L(\g p))_\mathrm{\bar l}$}
					\For{$k = 1,\dots,N_\mathrm d$}
						\If{$e_\mathrm{\bar \imath}M_\mathrm{\bar l, j(\bar \imath)}e_\mathrm k M_\mathrm{\bar l,j(k)} > 0$}
							\If{$\abs {e_k} > 0.01$} \Comment{Small errors can increase modifying the parameter}
								\State $\hat \nabla \mathcal L = \hat \nabla \mathcal L -|M_\mathrm{\bar l,j(k)}|$
							\EndIf
						\Else
							\State $\hat \nabla \mathcal L = \hat \nabla \mathcal L +|M_\mathrm{\bar l,j(k)}|$
						\EndIf
					\EndFor
					\If{$\hat \nabla \mathcal L < 0$}
						\State $p_\mathrm{\bar l} = \text{calibrate}(\bar l,p_\mathrm{\bar l},p^\mathrm R_\mathrm{\bar l},e_\mathrm{\bar \imath},M_\mathrm{\bar l, j(\bar \imath)})$ \Comment{Modify $p_\mathrm{\bar l}$}
						\State $it = it+1$
						\State Break
					\Else	\Comment{$e_\mathrm{\bar \imath}$ can not lower}
						\State $a_\mathrm{\bar \imath} = 0$
					\EndIf
				\EndWhile
			\EndIf
		\doWhile{$it \le it_\mathrm{max}$ \& $\mathcal L(\g p)\ge tol$}
	\end{algorithmic}
\end{algorithm}

Observe that if the sign of $e_\mathrm k M_\mathrm{\bar l,j(k)}$ is positive or negative, then $p_\mathrm{\bar l}$ has to be increased or decreased, respectively, to reduce the error on data $d_\mathrm{k}$. Therefore, if the sign of $e_\mathrm k M_\mathrm{\bar l,j(k)}$ is concordant to the sign of $e_\mathrm{\bar \imath}M_\mathrm{\bar l, j(\bar \imath)}$, modifying $p_\mathrm{\bar l}$ affects the errors $e_\mathrm{k}$ and $e_\mathrm{\bar \imath}$ in the same way.

The \emph{calibrate} function modifies the parameter $p_\mathrm{\bar l}$ sampling from a uniform distribution between the value of $p_\mathrm{\bar l}$ and its maximal or minimal bound (Appendix \ref{app:bounds}) if $e_\mathrm{\bar \imath}M_\mathrm{\bar l, j(\bar \imath)}$ is positive or negative, respectively. The bounds of the parameter depend on $p^\mathrm{R}_\mathrm{\bar l}$ and  the type of parameter considered (e.g., $\bar l$ can refer to an active elastance). The randomness of the step is due to the lack of the knowledge of an optimal step and the convergence of this calibration method can depend on it.

CMC offers an advantageous trade-off between accuracy and complexity. Its gradient-free nature allows to not compute the gradient at each iteration of the calibration procedure, thus reducing the computational cost. Nonetheless, performing a random step between two consecutive iterations, CMC could ignore local information worsening the convergence.

\subsection{Hybrid CMC-L-BFGS-B method} \label{sec:CMC-LBFGS_B}

The hybryd method CMC-L-BFGS-B, firstly, applies CMC to avoid local minima related to high values of the loss function and to get close to a better minimum; secondly, it applies L-BFGS-B to improve the estimate of this better minimum.

\section{Results}\label{sec:results}

In this section we report the results of the following tests for CMC, L-BFGS-B and CMC-L-BFGS-B:
\begin{enumerate}
\item robustness with respect to in silico generated data (Section \ref{sec:robCD});
\item robustness with respect to the initial guess of parameters (Section \ref{sec:robIG});
\item robustness with respect to noisy data (Section \ref{sec:robND});
\item calibration based on patient-specific data (Section \ref{sec:RP}).
\end{enumerate}
We run the first three tests on datasets of in silico generated data. We generate the data sampling the parameters in the associated ranges for the calibration procedure from a uniform distribution (Appendix \ref{app:bounds}) and computing the related model outputs.

\begin{figure}[H]
\begin{center}
\includegraphics[width=\linewidth]{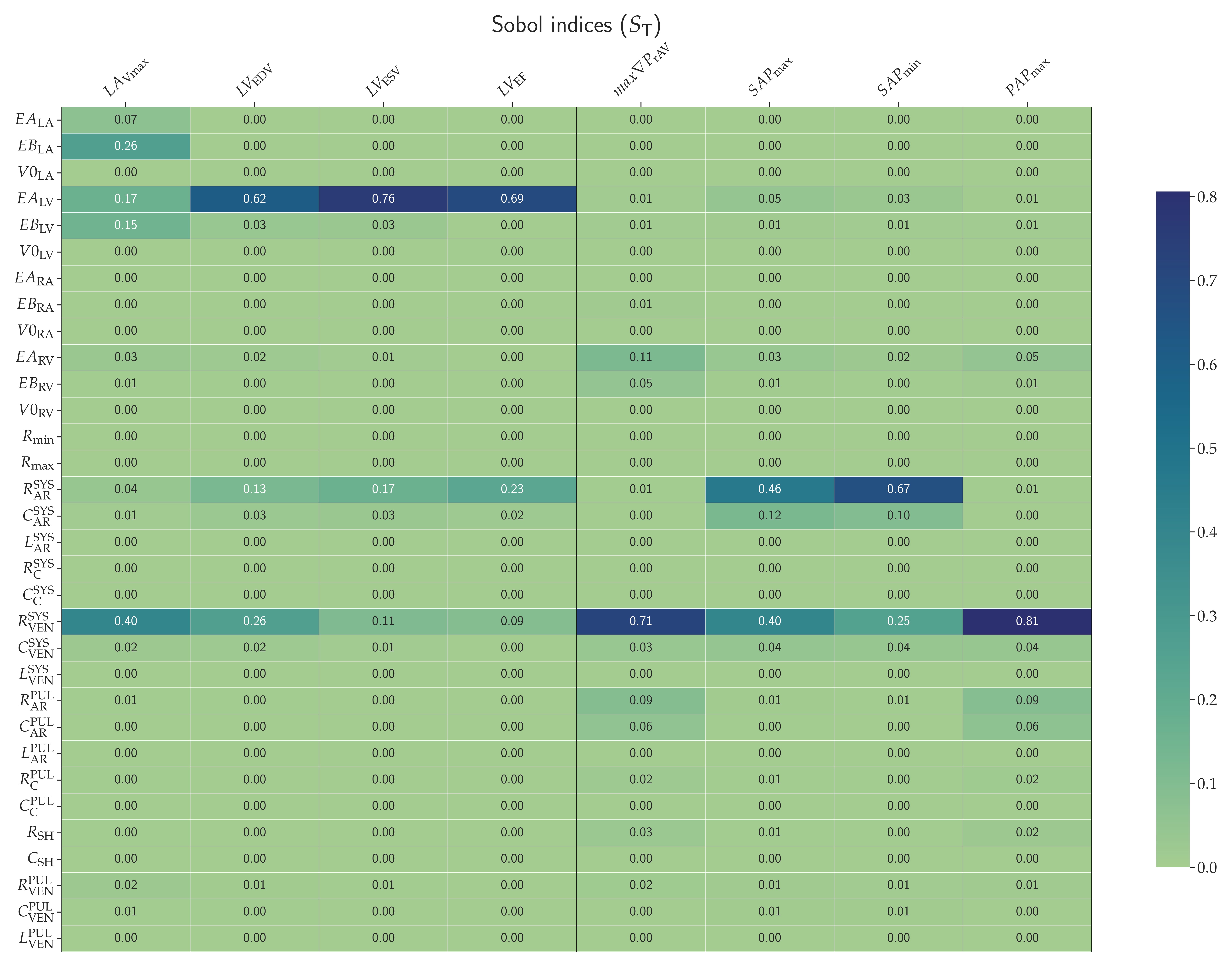}
\caption{Total-effect Sobol indices between parameters and model outputs related to data. A detailed definition of parameters and model outputs is provided in Table \ref{table:params} and Table \ref{table:qoi}.}
\label{fig:sobol}
\end{center}
\end{figure}

According to the data at disposal, we select for the calibration procedures the parameters related at least to one total Sobol index equal to or greater than $0.1$ (Figure \ref{fig:sobol}), namely: $EB_\mathrm{LA}$, $EA_\mathrm{LV}$, $EB_\mathrm{LV}$, $EA_\mathrm{RV}$, $R_\mathrm{AR}^\mathrm{SYS}$, $C_\mathrm{AR}^{SYS}$ and $R_\mathrm{VEN}^\mathrm{SYS}$ for a calibration procedure with a full set of data.

We implement CMC and L-BFGS-B in Matlab and in Python, respectively. We compute the numerical solution of the lumped-parameter cardiocirculatory model, in Matlab, by means of \emph{ode15s} solver which is a variable step, variable order solver, whereas, in Python, by means of Dormand-Prince method, which is an adaptive stepsize Runge-Kutta method. For each setting of the parameters here considered, we run the model for $25$ cycles until convergence to the regime solution of the dynamical system and we consider only the last heartbeat. To apply L-BFGS-B, we compute the exact gradient of the loss function by means of automatic differentiation. Since we are considering a loss function $\mathcal L: \R^\mathrm{N_\mathrm{p}}\rightarrow \R$ where $N_\mathrm{p}\ge 1$, it is always more convenient to use reverse (rather than forward) automatic differentiation, as it allows to compute the gradient in only one iteration instead of $N_\mathrm{p}$. To compute the reverse mode gradient, we use the Python library \emph{Jax} \cite{bradbury2018jax,salvador2023fast}. Once the gradient is at our disposal, we apply the quasi-Newton method L-BFGS-B to minimize the loss function.

We consider a calibration procedure successful if the root mean squared error 
\begin{align}
RMSE(\mathbf p) = \sqrt{MSE(\mathbf p)}
\end{align}
gets to a value lower than $10^{-1}$. Even if this condition is satisfied, the calibration method does not stop to possibly obtain a better estimate of the parameters. CMC-L-BFGS-B performs some steps of CMC to obtain a $MSE$ lower than $2.5\cdot 10^{-2}$ and then it applies L-BFGS-B. 

\subsection{Robustness and accuracy on silico generated data} \label{sec:robCD}

In this section, we test the robustness of the calibration methods in estimating the parameters on a dataset of $20$ different in silico generated data. The initial setting of parameters for the calibration procedures is displayed in Table \ref{table:params}.

\begin{figure}[H]
	\centering
 	\subfigure[]{\includegraphics[width=0.49\linewidth]{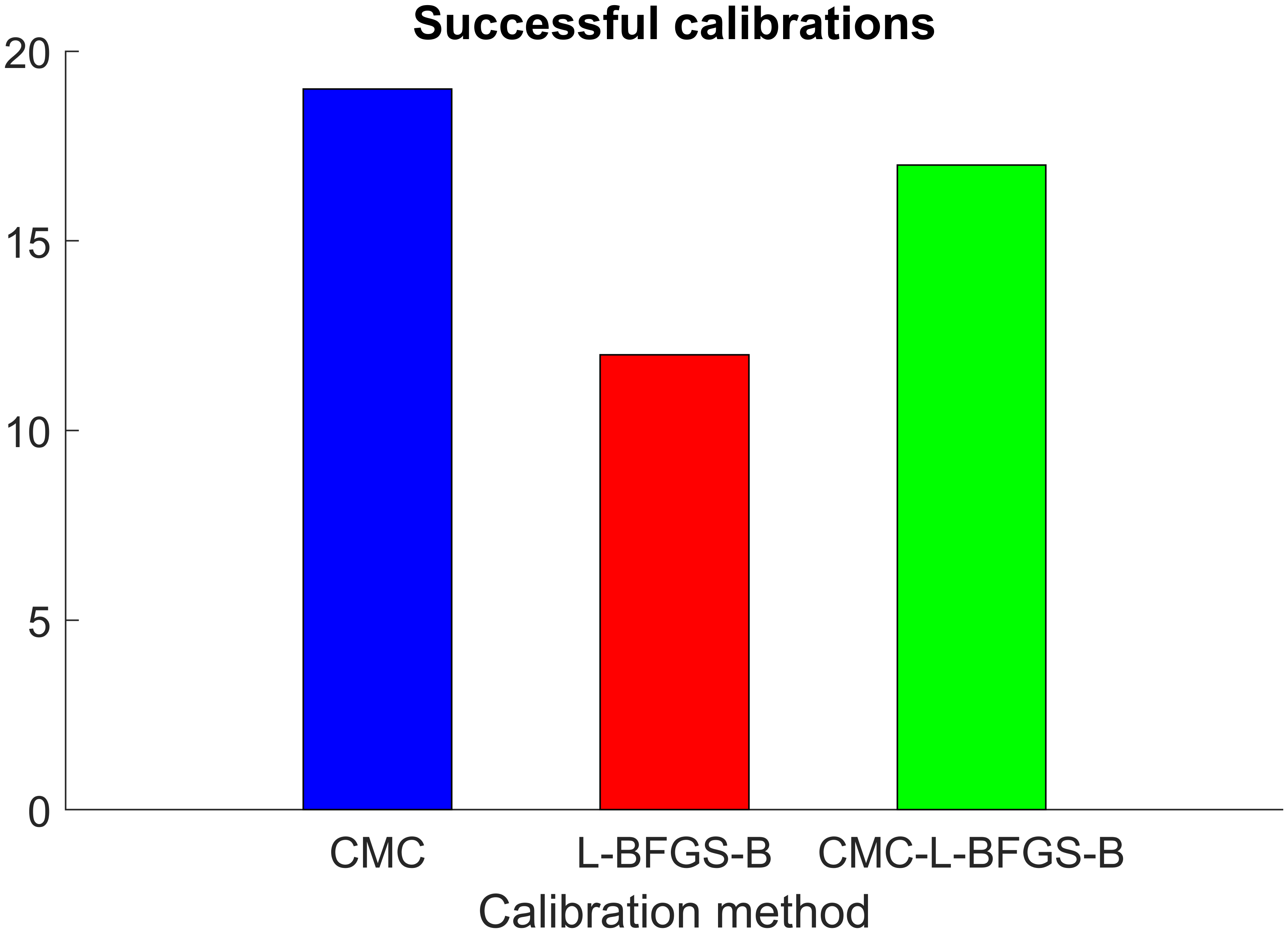}} 
	\subfigure[]{\includegraphics[width=0.49\linewidth]{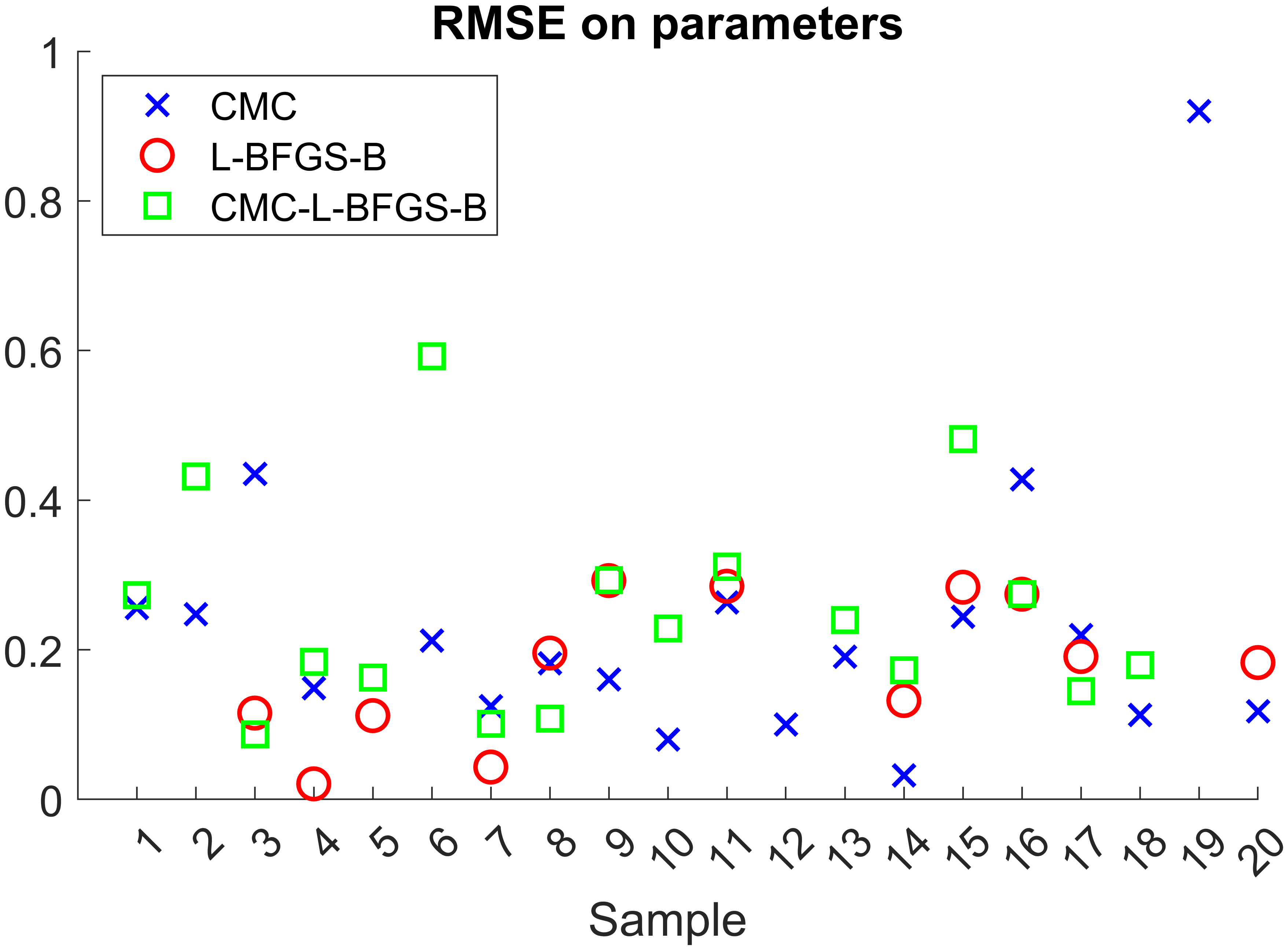}}
	\caption{Test 1. Number of successful calibration procedures for each calibration method (a) and $RMSE$ between estimated and real parameters for each sample (b). Only succesful calibrations are reported in the figure.}
	\label{fig:robCD}
\end{figure}

CMC is successful for $19$ samples, whereas L-BFGS-B is successful for $12$ samples and CMC-L-BFGS-B for $17$ samples (Figure \ref{fig:robCD}.a), assuming an intermediate behaviour between the two original calibration methods. Therefore, the CMC step increases the number of successful calibration procedures. Despite this advantage, L-BFGS-B estimates the parameters more precisely than the other two calibration methods obtaining on average an error that is lower than the ones returned by the other two calibration methods (Figure \ref{fig:robCD}.b). CMC returns an $RMSE$ on the parameters that is higher than the one returned by L-BFGS-B, highlighting that CMC is charachterized by a higher number of successful calibration procedures at the cost of a worse parameter estimation.

Finally, we measured the computational times of the three calibration methods on a standard laptop (AMD Ryzen 7 2700U, $2.20$GHz, $16$GB RAM):
\begin{itemize}
\item the computation of the correlation matrix takes $2.2$h and CMC takes on average $7.5$min for sample, for a total time of $4.7$h;
\item L-BFGS-B takes on average $57$min for sample, for a total time of $19$h;
\item CMC-L-BFGS-B takes on average $65$min for sample, for a total time, considering the correlation matrix computation, of $24$h split in the following way: $2.2$h for the computation of the correlation matrix; $9$ min for the CMC step; $22$h for the L-BFGS-B step.
\end{itemize}

CMC performs $4$ times better than L-BFGS-B in terms of computational time. Moreover, not considering the computation of the correlation matrix of CMC, it performs $7.5$ times better than L-BFGS-B.

\subsection{Robustness with respect to the initial guess of parameters} \label{sec:robIG}

In this section, we test the robustness of the calibration methods on sample $7$ of the previous Section with respect to the initial guess of parameters. We start the calibration procedures from $19$ different randomly selected initial settings of parameters lying in the associated ranges for the calibration procedure. Sample $7$ was associated with a successful calibration procedure for each calibration method in the previous Section. We want to observe if the three calibration methods converge to the same setting of parameters or to different ones.

\begin{figure}[H]
	\centering
 	\subfigure[]{\includegraphics[width=0.47\linewidth]{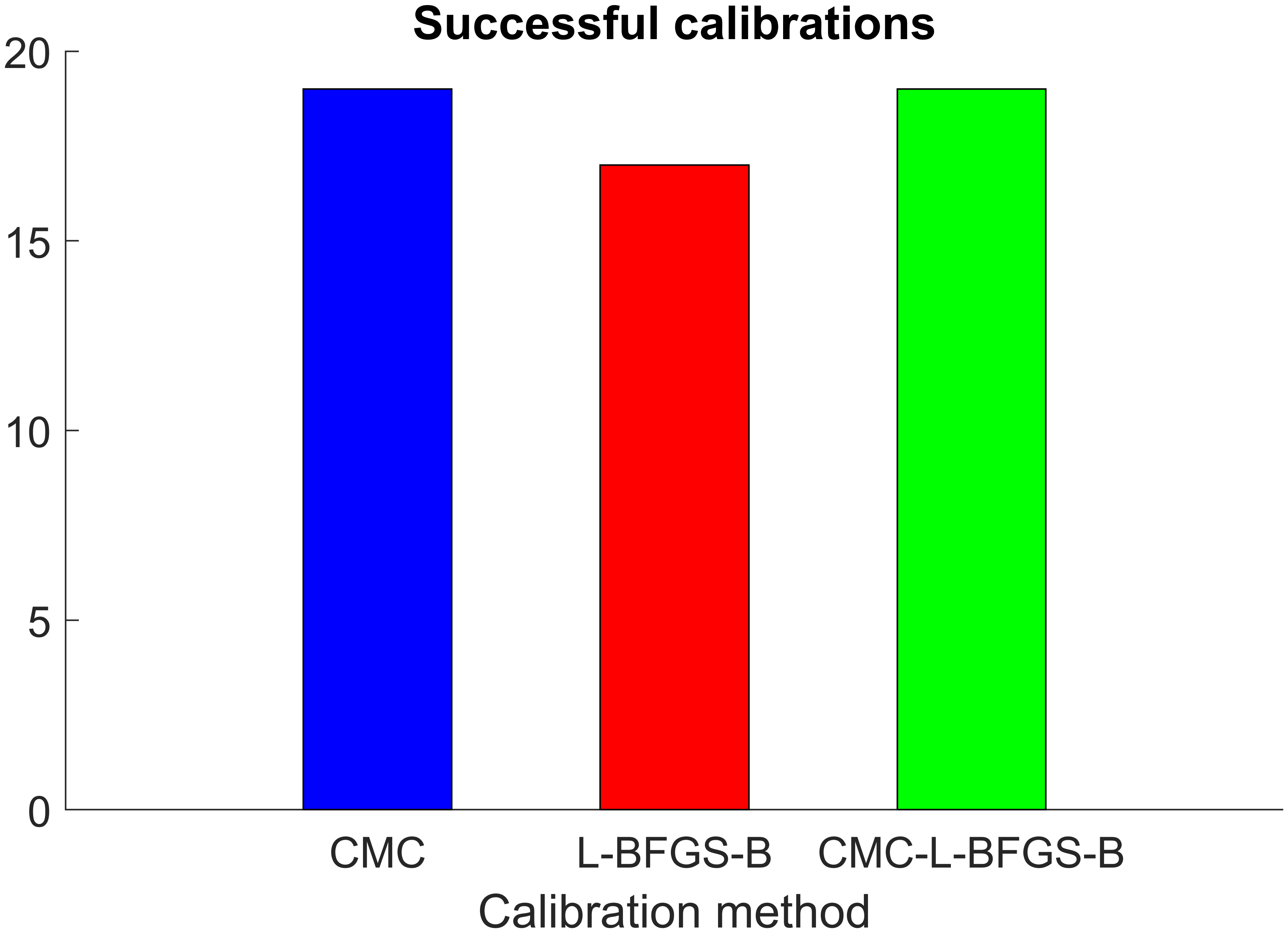}} 
	\subfigure[]{\includegraphics[width=0.47\linewidth]{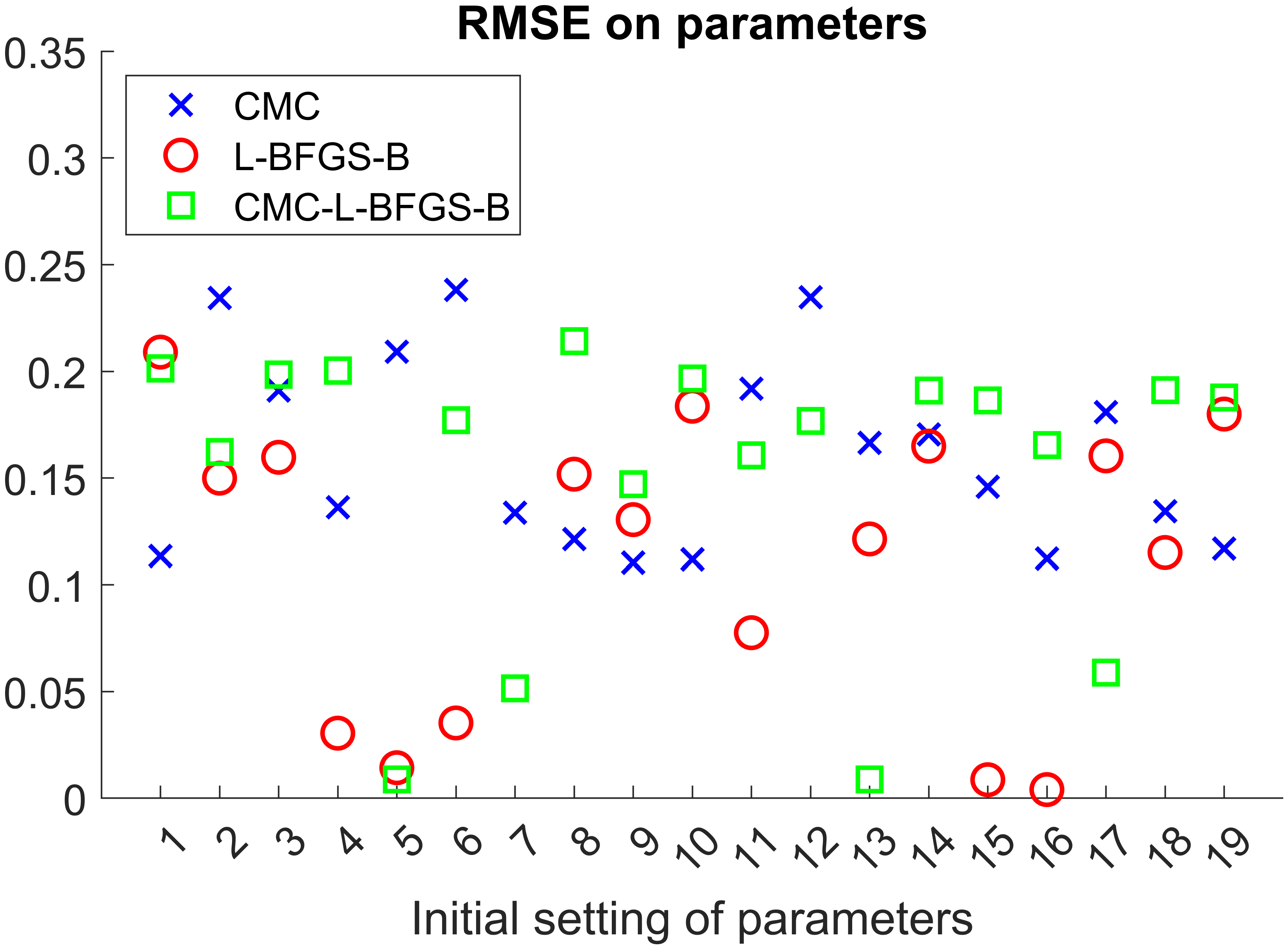}}
 	\subfigure[]{\includegraphics[width=0.47\linewidth]{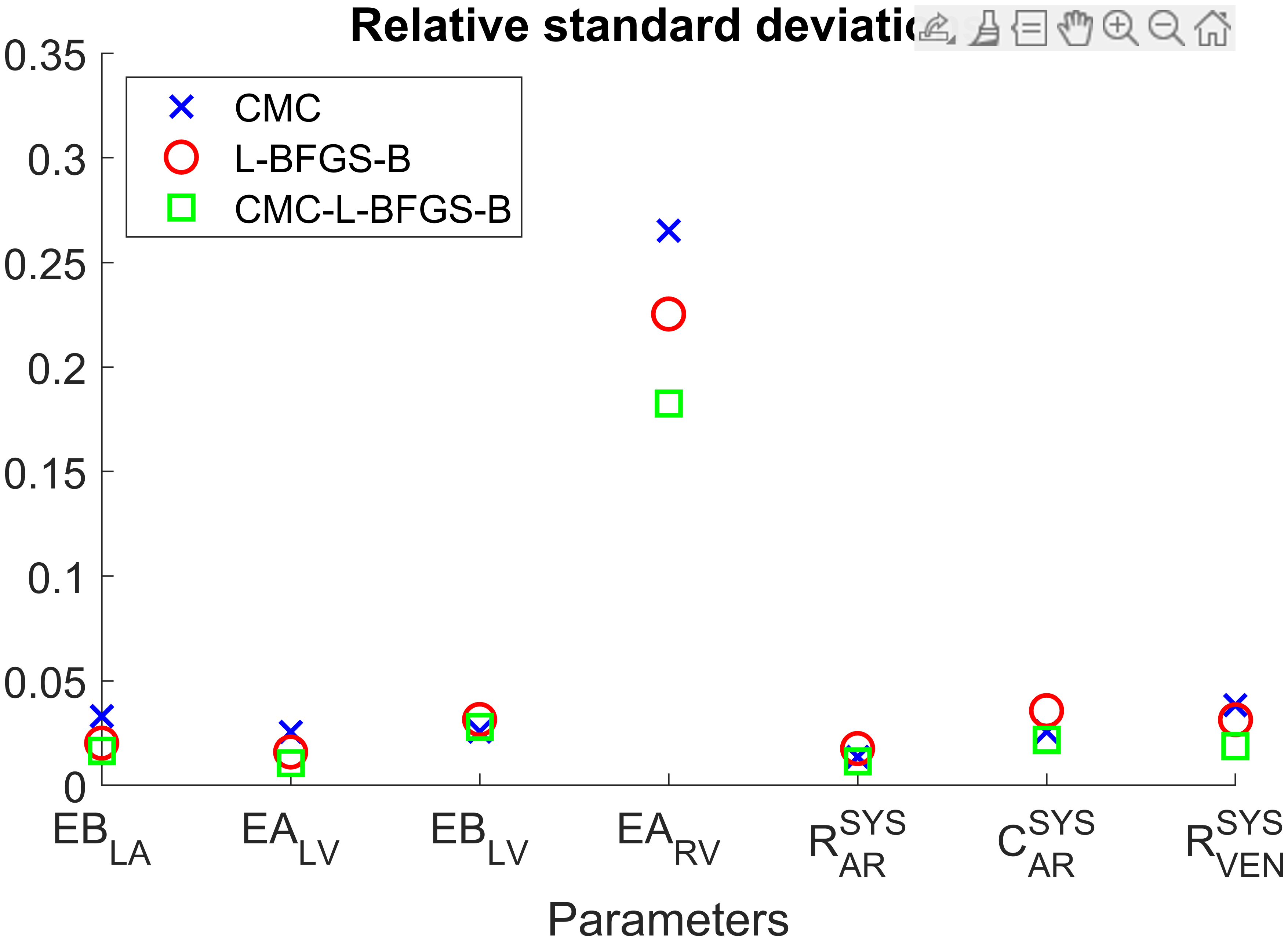}}
	\caption{Test 2. Number of successful calibration procedures for each calibration method (a), $RMSE$ between estimated and real parameters for each initial guess of parameters (b) and relative standard deviations of the estimated parameters (c). Only succesful calibrations are reported in the figure.}
	\label{fig:robIG}
\end{figure}

As preliminary results, the calibration procedure of the model with CMC is successful for $19$ settings, whereas L-BFGS-B is successful for $17$ settings and CMC-L-BFGS-B for $19$ settings (Figure \ref{fig:robIG}.a). With respect to the previous section, the differences in the parameter errors between each calibration method are smaller (Figure \ref{fig:robIG}.b).

We compute for each parameter the relative standard deviation of the estimated parameters with respect to the real value of the parameter (Figure \ref{fig:robIG}.c) to determine the robustness of the calibration methods with respect to the initial guess. Each calibration method achieves a relative standard deviation lower than $4\%$ for the estimated parameters different from $EA_\mathrm{RV}$. For what concerns $EA_\mathrm{RV}$:
\begin{itemize}
\item CMC achieves a relative standard deviation of $27\%$;
\item L-BFGS-B achieves a relative standard deviation of $23\%$;
\item CMC-L-BFGS-B achieves a relative standard deviation of $18\%$.
\end{itemize}

The relative standard deviation of the estimate of $EA_\mathrm{RV}$ is larger than the ones of the other parameters because $EA_\mathrm{RV}$ affects significantly only $\max \nabla p_\mathrm{rAV}$ (total Sobol index greater than $0.1$ in Figure \ref{fig:sobol}) and the related Sobol index is small ($0.11$), limiting the information to determine it accurately. CMC-L-BFGS-B achieves the best results in terms of relative standard deviations, even if the $RMSE$ on the parameters are mostly higher than the ones returned by L-BFGS-B (Figure \ref{fig:robIG}.b). Nonetheless, except for $EA_\mathrm{RV}$, all the three calibration methods are robust with respect to the initial guess of parameters.

\subsection{Robustness with respect to noisy data} \label{sec:robND}
In this section we apply syntethic noise to the sample $7$ of Section \ref{sec:robCD}, building a dataset of $20$ samples. Noisy data represent a realistic setting and the capacity to correctly estimate the parameters in this situation is crucial for calibration methods. We repeat the calibration procedure for each sample to test the robustness of the calibration methods with respect to noisy data. We sample the noise of each data from a normal distribution with zero mean and standard deviation equal to the expected measurement error associated to the data (Table \ref{table:measErr}). The initial setting of parameters for the calibration procedures is the same of Section \ref{sec:robCD}.

\begin{figure}[H]
	\centering
 	\subfigure[]{\includegraphics[width=0.49\linewidth]{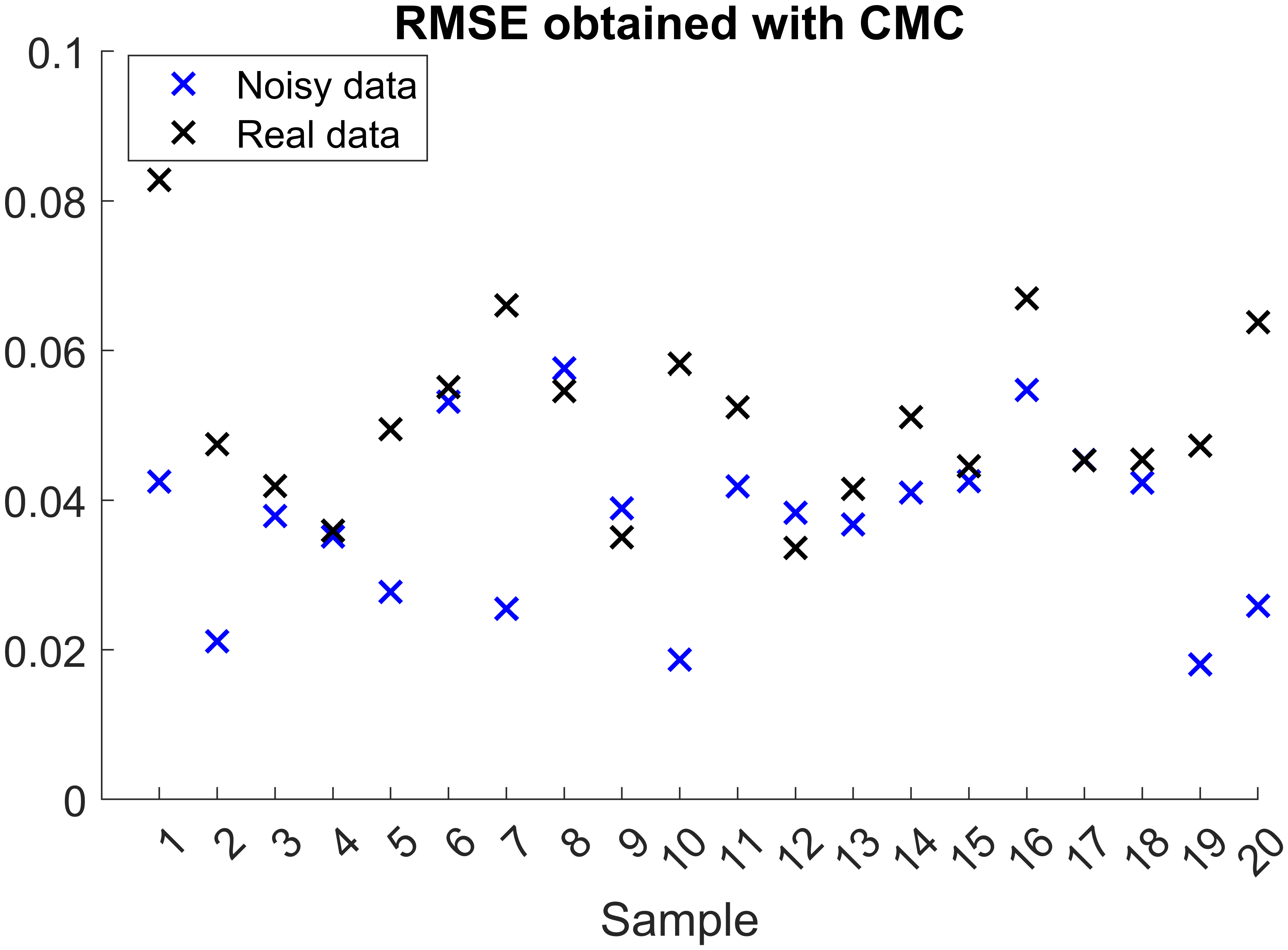}} 
	\subfigure[]{\includegraphics[width=0.49\linewidth]{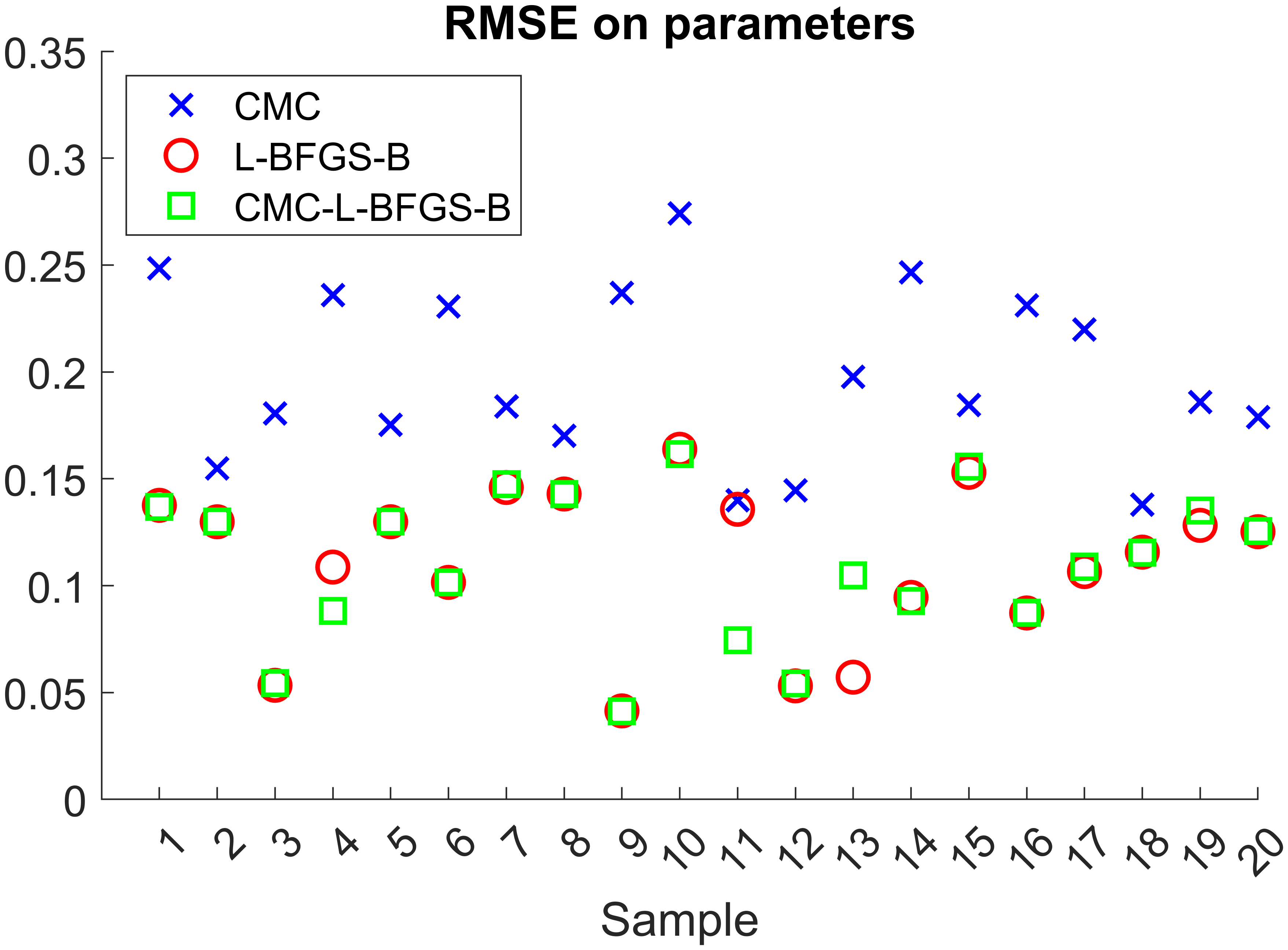}}
	\caption{Test 3. $RMSE$ on noisy and actual in silico generated data for CMC (a) and relative standard deviations of the estimated parameters for the $20$ samples (b).}
	\label{fig:robN}
\end{figure}

For each of the three calibration methods, all the $20$ calibration procedures are successful. Computing the $RMSE$ between model outputs and actual data (unaffected by noise), we obtain, for almost all the samples, higher errors than the $RMSE$ computed with respect to noisy data. Nevertheless, the $RMSE$ with respect to actual data is less than $10^{-1}$, i.e. the threshold for a successful calibration procedure (Figure \ref{fig:robN}.a in the case of CMC).

The relative standard deviations (Figure \ref{fig:robN}.b) of each parameter with respect to its value are lower than $20\%$ except for the one related to $EA_\mathrm{RV}$ and CMC. The relative standard deviations are, in general, higher with respect to the ones in Figure \ref{fig:robIG}.c because in this case the calibration procedures aim to converge to noisy data and not to the real ones. So, the parameters that achieve the minimum of the noisy loss function are not the real ones. Therefore, a calibration method less accurate than another one in estimating the parameters returns a worse estimate of noisy parameters, but it could achieve a better result on the actual ones. Indeed, CMC achieves the best results in terms of relative standard deviations for $5$ out of $7$ parameters despite being the least accurate calibration method in estimating the noisy parameters among the three calibration methods, as seen in Section \ref{sec:robCD}.

\subsection{Calibrations on clinical data}\label{sec:RP}

In this section, we perform the calibration procedures of the model on clinical data of two patients affected by COVID-19 related pneumonia (Table \ref{table:real_test}). Centro Cardiologico Monzino in Milan and L. Sacco Hospital provided the data. 
The patient from Monzino has a restricted number of clinical data, so we choose the parameters to calibrate accordingly to the global sensitivity analysis (as described at the beginning of Section \ref{sec:results}): namely, $EA_\mathrm{LV}$, $EA_\mathrm{RV}$, $R_\mathrm{AR}^\mathrm{SYS}$, $C_\mathrm{AR}^\mathrm{SYS}$ and $R_\mathrm{VEN}^\mathrm{SYS}$. We compare the estimated parameters and the ventricular PV (Pressure-Volume) loops that give clinical insight of the cardiac condition of a patient. 

All the loss functions are less than $10^{-1}$ (Table \ref{table:losses}), so each calibration procedure is successful. For both patients CMC-L-BFGS-B and L-BFGS-B obtain the best results in terms of loss function. 

\begin{figure}[H]
	\centering
 	\subfigure[]{\includegraphics[width=0.49\linewidth]{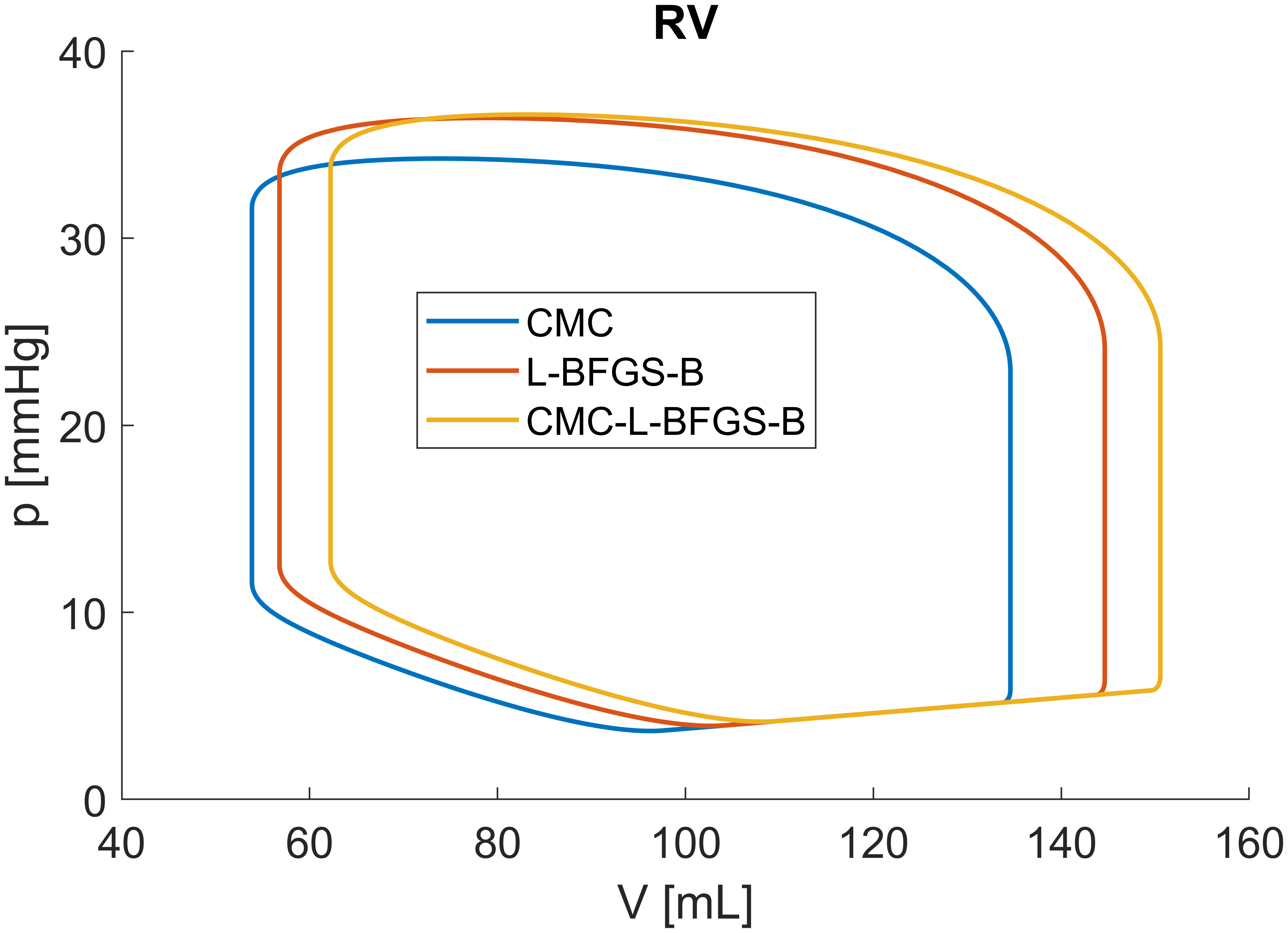}} 
	\subfigure[]{\includegraphics[width=0.49\linewidth]{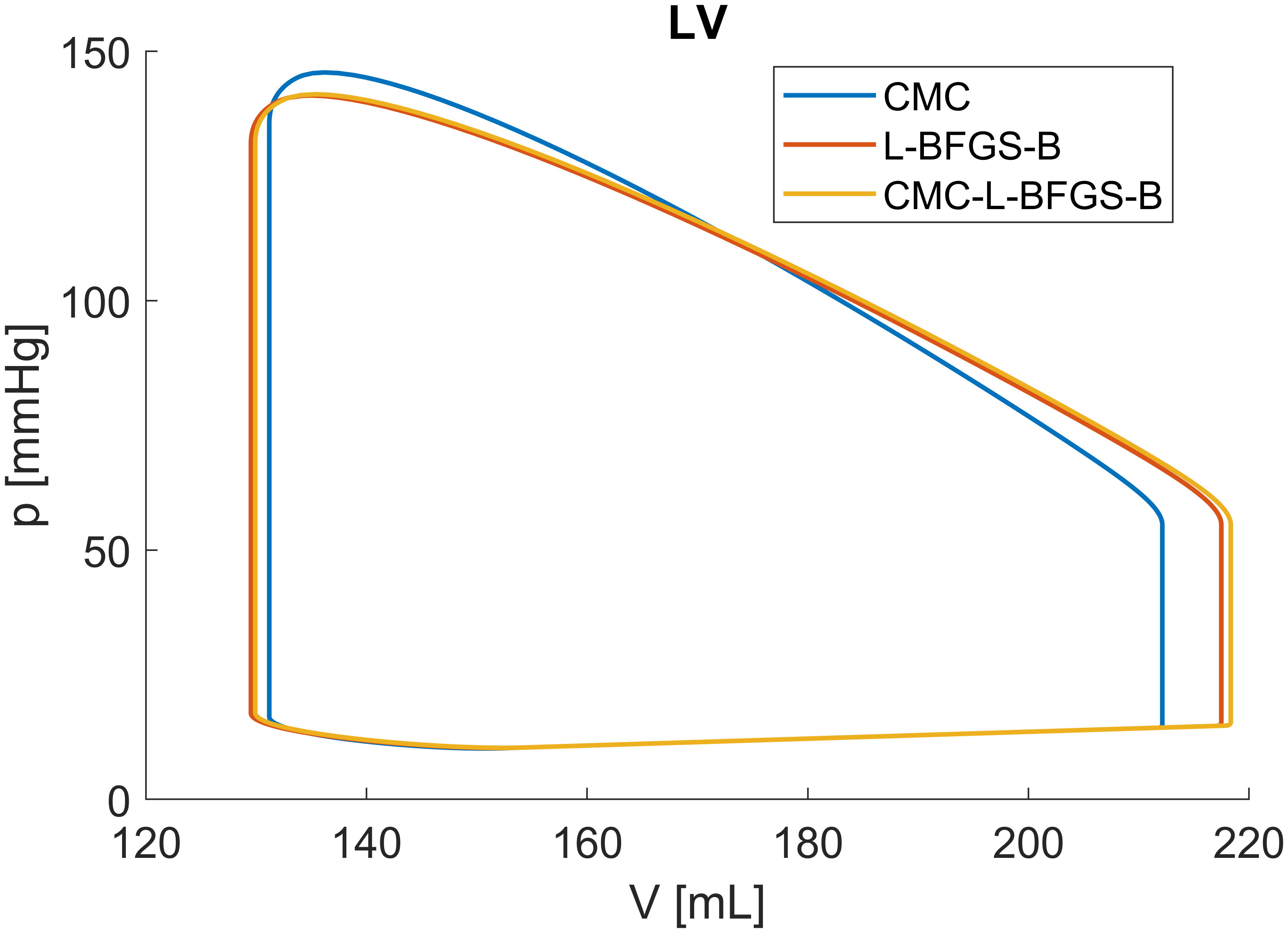}}
	
 	\subfigure[]{\includegraphics[width=0.49\linewidth]{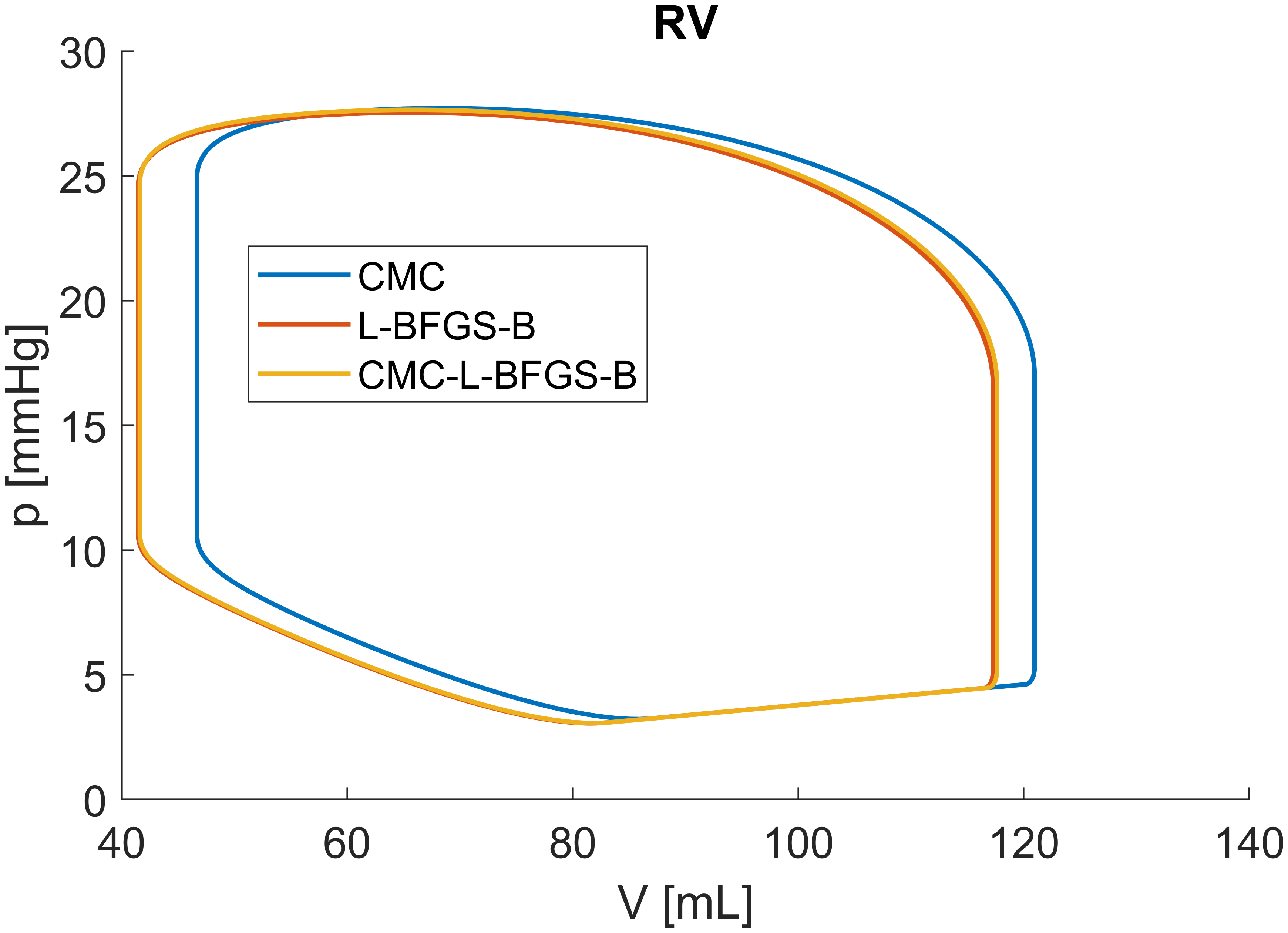}} 
	\subfigure[]{\includegraphics[width=0.49\linewidth]{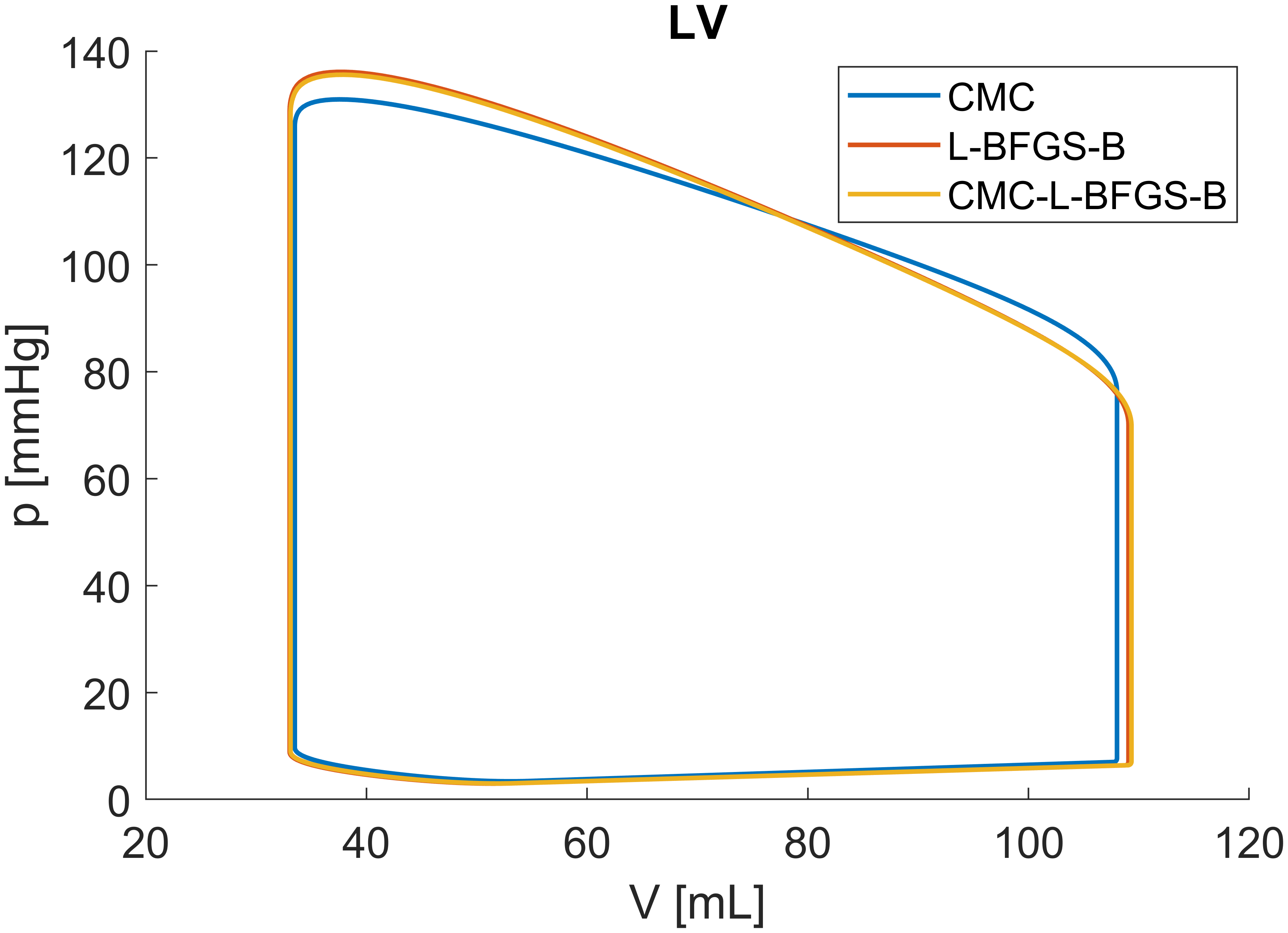}}
	\caption{Test 4. Estimated right ventricular (a) and left ventricular (b) PV loops of the Monzino patient and right ventricular (c) and left ventricular (d) PV loops of the Sacco patient.}
	\label{fig:PVloops}
\end{figure}

In the case of the patient from Monzino hospital, the relative standard deviations between the three calibration methods of the estimated parameters are less than $10\%$ (Table \ref{table:est_params}), indicating similar results for the three calibration methods. L-BFGS-B and CMC-L-BFGS-B estimate the same values of parameters except for $EA_\mathrm{RV}$. Therefore, the left ventricular PV loops are nearly the same for these two calibration methods, whereas for the right ventricular PV loops the differences are bigger (Figures \ref{fig:PVloops}.a and \ref{fig:PVloops}.b). CMC returns different estimates of the parameters related to the blood vessels with respect to the other two calibration methods. As a consequence, the PV loops returned from the former calibration method are slightly different from the other two. Nonetheless, from a clinical point of view, these differences are negligible. In particular, even if we only know the value of max $\nabla P_\mathrm{rAV}$ that is related to the right ventricle, the three calibration methods return the same condition: the $RV_{\mathrm{Pmax}}$ values (all greater than $30$ mmHg) are higher than the associated healthy range \cite{hurst1990heart}. This indicates an increased workload of the right ventricle. Instead, the values of $RV_{\mathrm{Pmin}}$ (all less than $8$ mmHg) are in healthy ranges (Table \ref{table:qoi}). For what concerns the left ventricular PV loops, the values of $LV_\mathrm{ESV}$ (all greater than $52$ mL), $LV_\mathrm{Pmax}$ (all greater than $140$ mmHg) and $LV_\mathrm{Pmax}$ (all greater than $9$ mmHg) are higher than the associated healthy range (Table \ref{table:qoi}), indicating an impaired left ventricular function. We do not compare the values of $RV_\mathrm{EDV}$, $RV_\mathrm{ESV}$ and $LV_\mathrm{EDV}$ with their indexed healthy ranges (Table \ref{table:qoi}) because we do not know the BSA of the patient at hand.

In the case of the patient from L. Sacco Hospital, the relative standard deviations of the estimated parameters are less than $12\%$ (Table \ref{table:est_params}). The discrepancies among the PV loops are diminished with respect to the previous case (Figure \ref{fig:PVloops}.c and \ref{fig:PVloops}.d), because more clinical data are available. The three right and left ventricular PV loops indicate that the patient, even if it is affected by COVID-19 related pneumonia, is in healthy conditions.

\section{Conclusions} \label{sec:concl}

This study introduces a new gradient-free calibration method for a cardiocirculatory mathematical model based on the correlation matrix between parameters and model outputs. 

We showed that the new method (named CMC) outperforms L-BFGS-B in terms of the number of successful calibrations (Section \ref{sec:robCD}). Thanks to its random step procedure, CMC avoids getting stuck in local minima and explores better than L-BFGS-B the hyperbox where the parameters can vary. This property echoes on the combined CMC-L-BFGS-B that also outperforms L-BFGS-B in terms of the number of successful calibrations. Because of its gradient-free nature, CMC returns a worse estimate of the parameters (Figure \ref{fig:robCD}.b). Nonetheless, we showed that, even if there are differences in the estimated parameters with the three calibration methods, the PV loops obtained for real patients indicate the same cardiac condition (Section \ref{sec:RP}).

Moreover, we showed that CMC outperforms both L-BFGS-B and CMC-L-BFGS-B in terms of computational time, offering a good trade-off between accuracy in estimating clinical data and computational cost (Section \ref{sec:robCD}). CMC-L-BFGS-B is the most robust method among the three with respect to the initial guess of parameters (Section \ref{sec:robIG}), even if it is the slowest due to the L-BFGS-B step.

Finally, CMC is the most robust of the three methods with respect to noisy data (Section \ref{sec:robND}). Due to the less accurate estimate of the parameters (Section \ref{sec:robCD}), CMC returns parameters related to the noisy data worse than the other two calibration methods. This enhances the estimate of the real parameters and makes CMC an effective choice in presence of noisy data, a typical situation that arises with clinical data.

We now discuss the limitations of CMC. First, the correlation coefficient could overlook strong non linear relationships between parameters and model outputs, thus neglecting some dependencies. We faced this problem, implementing CMC-L-BFGS-B to account for non linear relationships with L-BFGS-B. Nonetheless, we showed that CMC-L-BFGS-B features a lower convergence rate with respect to the other two calibration methods. Tuning the threshold of the $MSE$ (set to $2.5\cdot 10^\mathrm{-2}$ in this work) employed to switch CMC with L-BFGS-B could improve the performances of CMC-L-BFGS-B.

Secondly, the random step of CMC can both improve and worsen the value of the loss function avoiding not only local minima, but also minima associated with successful calibrations. An enhanced choice of the random step, for example reducing its maximum value according to the value of the loss function in the current parameter setting, could improve the convergence of CMC.

\bibliographystyle{unsrt}
\bibliography{Tonini_et_al_archive}

\begin{thebibliography}{10}

\bibitem{shi2011review}
Yubing Shi, Patricia Lawford, and Rodney Hose.
\newblock Review of zero-d and 1-d models of blood flow in the cardiovascular
  system.
\newblock {\em Biomedical engineering online}, 10:1--38, 2011.

\bibitem{quarteroni2016geometric}
Alfio Quarteroni, Alessandro Veneziani, and Christian Vergara.
\newblock Geometric multiscale modeling of the cardiovascular system, between
  theory and practice.
\newblock {\em Computer Methods in Applied Mechanics and Engineering},
  302:193--252, 2016.

\bibitem{dede2021modeling}
Luca Ded{\`e}, Francesco Regazzoni, Christian Vergara, Paolo Zunino, Marco
  Guglielmo, Roberto Scrofani, Laura Fusini, Chiara Cogliati, Gianluca Pontone,
  and Alfio Quarteroni.
\newblock Modeling the cardiac response to hemodynamic changes associated with
  covid-19: a computational study.
\newblock {\em Mathematical Biosciences and Engineering}, 18(4):3364--3383,
  2021.

\bibitem{fernandes2021integrated}
Luciano~G Fernandes, Paulo~R Trenhago, Ra{\'u}l~A Feij{\'o}o, and Pablo~J
  Blanco.
\newblock Integrated cardiorespiratory system model with short timescale
  control mechanisms.
\newblock {\em International Journal for Numerical Methods in Biomedical
  Engineering}, 37(11):e3332, 2021.

\bibitem{de2006modelling}
Claudio De~Lazzari, Marek Darowski, Gianfranco Ferrari, Domenico~M Pisanelli,
  and G~Tosti.
\newblock Modelling in the study of interaction of hemopump device and
  artificial ventilation.
\newblock {\em Computers in biology and medicine}, 36(11):1235--1251, 2006.

\bibitem{shi2006numerical}
Yubing Shi and Theodosios Korakianitis.
\newblock Numerical simulation of cardiovascular dynamics with left heart
  failure and in-series pulsatile ventricular assist device.
\newblock {\em Artificial organs}, 30(12):929--948, 2006.

\bibitem{tonini2024mathematical}
Andrea Tonini, Christian Vergara, Francesco Regazzoni, Luca Dede’, Roberto
  Scrofani, Chiara Cogliati, and Alfio Quarteroni.
\newblock A mathematical model to assess the effects of covid-19 on the
  cardiocirculatory system.
\newblock {\em Scientific Reports}, 14(1):8304, 2024.

\bibitem{zhu2023computational}
Ge~Zhu, Susree Modepalli, Mohan Anand, and He~Li.
\newblock Computational modeling of hypercoagulability in covid-19.
\newblock {\em Computer Methods in Biomechanics and Biomedical Engineering},
  26(3):338--349, 2023.

\bibitem{herrmann2020modeling}
Jacob Herrmann, Vitor Mori, Jason~HT Bates, and B{\'e}la Suki.
\newblock Modeling lung perfusion abnormalities to explain early covid-19
  hypoxemia.
\newblock {\em Nature communications}, 11(1):4883, 2020.

\bibitem{laubscher2023estimation}
Ryno Laubscher, Johan Van Der~Merwe, Philip Herbst, and Jacques Liebenberg.
\newblock Estimation of simulated left ventricle elastance using lumped
  parameter modelling and gradient-based optimization with forward-mode
  automatic differentiation based on synthetically generated noninvasive data.
\newblock {\em Journal of Biomechanical Engineering}, 145(2):021008, 2023.

\bibitem{bjordalsbakke2022parameter}
Nikolai~L Bj{\o}rdalsbakke, Jacob~T Sturdy, David~R Hose, and Leif~R Hellevik.
\newblock Parameter estimation for closed-loop lumped parameter models of the
  systemic circulation using synthetic data.
\newblock {\em Mathematical biosciences}, 343:108731, 2022.

\bibitem{saxton2023personalised}
Harry Saxton, Torsten Schenkel, Ian Halliday, and Xu~Xu.
\newblock Personalised parameter estimation of the cardiovascular system:
  Leveraging data assimilation and sensitivity analysis.
\newblock {\em Journal of Computational Science}, 74:102158, 2023.

\bibitem{saltelli2002making}
Andrea Saltelli.
\newblock Making best use of model evaluations to compute sensitivity indices.
\newblock {\em Computer physics communications}, 145(2):280--297, 2002.

\bibitem{sobol1993sensitivity}
IM~Sobo{\'l}.
\newblock Sensitivity estimates for nonlinear mathematical models.
\newblock {\em Math. Model. Comput. Exp.}, 1:407, 1993.

\bibitem{byrd1995limited}
Richard~H Byrd, Peihuang Lu, Jorge Nocedal, and Ciyou Zhu.
\newblock A limited memory algorithm for bound constrained optimization.
\newblock {\em SIAM Journal on scientific computing}, 16(5):1190--1208, 1995.

\bibitem{zhu1997algorithm}
Ciyou Zhu, Richard~H Byrd, Peihuang Lu, and Jorge Nocedal.
\newblock Algorithm 778: L-bfgs-b: Fortran subroutines for large-scale
  bound-constrained optimization.
\newblock {\em ACM Transactions on mathematical software (TOMS)},
  23(4):550--560, 1997.

\bibitem{dandel2022heart}
Michael Dandel.
\newblock Heart--lung interactions in covid-19: prognostic impact and
  usefulness of bedside echocardiography for monitoring of the right ventricle
  involvement.
\newblock {\em Heart Failure Reviews}, 27(4):1325--1339, 2022.

\bibitem{park2020eye}
John~F Park, Somanshu Banerjee, and Soban Umar.
\newblock In the eye of the storm: the right ventricle in covid-19.
\newblock {\em Pulmonary Circulation}, 10(3):2045894020936660, 2020.

\bibitem{diaz2021association}
Carlos Diaz-Arocutipa, Jose Saucedo-Chinchay, and Edgar Argulian.
\newblock Association between right ventricular dysfunction and mortality in
  covid-19 patients: A systematic review and meta-analysis.
\newblock {\em Clinical Cardiology}, 44(10):1360--1370, 2021.

\bibitem{mauri2020potential}
Tommaso Mauri, Elena Spinelli, Eleonora Scotti, Giulia Colussi, Maria~Cristina
  Basile, Stefania Crotti, Daniela Tubiolo, Paola Tagliabue, Alberto Zanella,
  Giacomo Grasselli, et~al.
\newblock Potential for lung recruitment and ventilation-perfusion mismatch in
  patients with the acute respiratory distress syndrome from coronavirus
  disease 2019.
\newblock {\em Critical care medicine}, 48(8):1129--1134, 2020.

\bibitem{peng2020using}
Qian-Yi Peng, Xiao-Ting Wang, Li-Na Zhang, and Chinese Critical Care Ultrasound
  Study~Group (CCUSG).
\newblock Using echocardiography to guide the treatment of novel coronavirus
  pneumonia, 2020.

\bibitem{regazzoni2022cardiac}
Francesco Regazzoni, Matteo Salvador, Pasquale~Claudio Africa, Marco Fedele,
  Luca Ded{\`e}, and Alfio Quarteroni.
\newblock A cardiac electromechanical model coupled with a lumped-parameter
  model for closed-loop blood circulation.
\newblock {\em Journal of Computational Physics}, 457:111083, 2022.

\bibitem{velthuis2015pulmonary}
Sebastiaan Velthuis, Veronique~MM Vorselaars, Cornelis~JJ Westermann, Repke~J
  Snijder, Johannes~J Mager, and Martijn~C Post.
\newblock Pulmonary shunt fraction measurement compared to contrast
  echocardiography in hereditary haemorrhagic telangiectasia patients: time to
  abandon the 100\% oxygen method?
\newblock {\em Respiration}, 89(2):112--118, 2015.

\bibitem{saha2021correlation}
Biplab~K Saha, Sana Ghalib, Hau Chieng, Chad Pezzano, Darren Lydon, Paul
  Feustel, Thomas~C Smith, and Amit Chopra.
\newblock Correlation of respiratory physiologic parameters in mechanically
  ventilated coronavirus disease 2019 patients.
\newblock {\em Critical Care Explorations}, 3(1):e0328, 2021.

\bibitem{gattinoni2020covid}
Luciano Gattinoni, Silvia Coppola, Massimo Cressoni, Mattia Busana, Sandra
  Rossi, and Davide Chiumello.
\newblock Covid-19 does not lead to a “typical” acute respiratory distress
  syndrome.
\newblock {\em American journal of respiratory and critical care medicine},
  201(10):1299--1300, 2020.

\bibitem{liang2009multi}
Fuyou Liang, Shu Takagi, Ryutaro Himeno, and Hao Liu.
\newblock Multi-scale modeling of the human cardiovascular system with
  applications to aortic valvular and arterial stenoses.
\newblock {\em Medical \& biological engineering \& computing}, 47:743--755,
  2009.

\bibitem{albanese2016integrated}
Antonio Albanese, Limei Cheng, Mauro Ursino, and Nicolas~W Chbat.
\newblock An integrated mathematical model of the human cardiopulmonary system:
  model development.
\newblock {\em American Journal of Physiology-Heart and Circulatory
  Physiology}, 310(7):H899--H921, 2016.

\bibitem{zhang2009heart}
Gus~Q Zhang and Weiguo Zhang.
\newblock Heart rate, lifespan, and mortality risk.
\newblock {\em Ageing research reviews}, 8(1):52--60, 2009.

\bibitem{pearson1895vii}
Karl Pearson.
\newblock Vii. note on regression and inheritance in the case of two parents.
\newblock {\em proceedings of the royal society of London},
  58(347-352):240--242, 1895.

\bibitem{bradbury2018jax}
James Bradbury, Roy Frostig, Peter Hawkins, Matthew~James Johnson, Chris Leary,
  Dougal Maclaurin, George Necula, Adam Paszke, Jake VanderPlas, Skye
  Wanderman-Milne, et~al.
\newblock Jax: composable transformations of python+ numpy programs.
\newblock 2018.

\bibitem{salvador2023fast}
Matteo Salvador, Francesco Regazzoni, Alfio Quarteroni, et~al.
\newblock Fast and robust parameter estimation with uncertainty quantification
  for the cardiac function.
\newblock {\em Computer Methods and Programs in Biomedicine}, 231:107402, 2023.

\bibitem{hurst1990heart}
JW~Hurst, RC~Schlant, CE~Rackley, EH~Sonnenblick, and NK~Wenger.
\newblock The heart, arteries e veins, 1990.

\bibitem{lang2015recommendations}
Roberto~M Lang, Luigi~P Badano, Victor Mor-Avi, Jonathan Afilalo, Anderson
  Armstrong, Laura Ernande, Frank~A Flachskampf, Elyse Foster, Steven~A
  Goldstein, Tatiana Kuznetsova, et~al.
\newblock Recommendations for cardiac chamber quantification by
  echocardiography in adults: an update from the american society of
  echocardiography and the european association of cardiovascular imaging.
\newblock {\em European Heart Journal-Cardiovascular Imaging}, 16(3):233--271,
  2015.

\bibitem{rudski2010guidelines}
Lawrence~G Rudski, Wyman~W Lai, Jonathan Afilalo, Lanqi Hua, Mark~D
  Handschumacher, Krishnaswamy Chandrasekaran, Scott~D Solomon, Eric~K Louie,
  and Nelson~B Schiller.
\newblock Guidelines for the echocardiographic assessment of the right heart in
  adults: a report from the american society of echocardiography: endorsed by
  the european association of echocardiography, a registered branch of the
  european society of cardiology, and the canadian society of echocardiography.
\newblock {\em Journal of the American society of echocardiography},
  23(7):685--713, 2010.

\end{thebibliography}

\appendix

\section{Bounds on the parameters} \label{app:bounds}

In this Section we describe the bounds of the hyperbox used to perform the global sensitivity analysis and the calibration procedures of the lumped-parameter cardiocirculatory model. 

For what concerns the global sensitivity analysis, we have previously chosen a reference setting of parameters $\mathbf p^\mathrm{R}$ to model an ideal heatlhy individual (Section \ref{sec:model}) that differs from a generic healthy one. Therefore, the parameters can vary in a hyperbox to account for different healthy conditions: given the reference setting of parameters $\mathbf{p}^\mathrm{R}$ and fixing the times of the cardiac cycle according to the input $HR$, we sample the parameter $p_\mathrm{l}$ ($1\le l \le N_\mathrm{p}=32$) in $\left[\left(1-\frac{2}{3}\right)p^\mathrm{R}_\mathrm{l}, \left(1+\frac{2}{3}\right)p^\mathrm{R}_\mathrm{l}\right]$.

For what concerns the calibration methods, we fixed to the reference value the parameters related to total Sobol indices all less than $0.1$. The other parameters are free to vary in the same ranges as for the global sensitivity analysis. Moreover, if we suppose that the severe COVID-19-related pneumonia is ongoing, we will consider three additional changes:
\begin{itemize}
\item The active elastances of the four cardiac chambers can be further halved according to the impairment of the cardiac function due to the infection\cite{dandel2022heart}.
\item The pulmonary resistances can further triple according to the registered increase in pulmonary resistances reported in clinical literature \cite{mauri2020potential,peng2020using}. 
\item The pulmonary compliances can be further divided by $3$ according to the consequences of the endothelial damage of the pulmonary blood vessels\cite{mauri2020potential}.
\end{itemize}

\newpage

\pagenumbering{gobble}

\section*{Tables}
\thispagestyle{empty}
\begin{table}[H]
\begin{center}
\scriptsize
\begin{tabular}{|c|c|c|c|}
\hline
Parameter & Unit & Reference value & Description\\
\hline
$EA_{\mathrm{LA}}$ & mmHg/mL & $0.38$ & Left atrial active elastance\\
$EB_{\mathrm{LA}}$ & mmHg/mL & $0.27$ & Left atrial passive elastance\\
$V_{U,\mathrm{LA}}$ & mL & $2.31$ & Left atrial unloaded volume\\
$EA_{\mathrm{LV}}$ & mmHg/mL & $2.7$ & Left ventricular active elastance\\
$EB_{\mathrm{LV}}$ & mmHg/mL & $0.069$ & Left ventricular passive elastance\\
$V_{U,\mathrm{LV}}$ & mL & $3.54$ & Left ventricular unloaded volume\\
$EA_{\mathrm{RA}}$ & mmHg/mL & $0.13$ & Right atrial active elastance\\
$EB_{\mathrm{RA}}$ & mmHg/mL & $0.20$ & Right atrial passive elastance\\
$V_{U,\mathrm{RA}}$ & mL & $3.54$ & Right atrial unloaded volume\\
$EA_{\mathrm{RV}}$ & mmHg/mL & $0.43$ & Right ventricular active elastance\\
$EB_{\mathrm{RV}}$ & mmHg/mL & $0.041$ & Right ventricular passive elastance\\
$V_{U,\mathrm{RV}}$ & mL & $8.41$ & Right ventricular unloaded volume\\
$R_{\mathrm{min}}$ & mmHg$\cdot$ s/mL & 0.0063 & Minimal valve resistance\\
$R_{\mathrm{max}}$ & mmHg$\cdot$ s/mL & 94168 & Maximal valve resistance\\
$R_{\mathrm{AR}}^{\mathrm{SYS}}$ & mmHg$\cdot$ s/mL & 0.59 & Systemic arterial resistance\\
$C_{\mathrm{AR}}^{\mathrm{SYS}}$ & mL/mmHg & 1.33 & Systemic arterial compliance\\
$L_{\mathrm{AR}}^{\mathrm{SYS}}$ & mmHg$\cdot$ s$^\mathrm{2}$/mL & 0.00021 & Systemic arterial inertia\\
$R_{\mathrm{C}}^{\mathrm{SYS}}$ & mmHg$\cdot$ s/mL & 0.022 & Systemic capillary resistance\\
$C_{\mathrm{C}}^{\mathrm{SYS}}$ & mL/mmHg & 0.28 & Systemic capillary compliance\\
$R_{\mathrm{VEN}}^{\mathrm{SYS}}$ & mmHg$\cdot$ s/mL & 0.36 & Systemic venous resistance\\
$C_{\mathrm{VEN}}^{\mathrm{SYS}}$ & mL/mmHg & 75 & Systemic venous compliance\\
$L_{\mathrm{VEN}}^{\mathrm{SYS}}$ & mmHg$\cdot$ s$^\mathrm{2}$/mL & 0.000021 & Systemic venous inertia\\
$R_{\mathrm{AR}}^{\mathrm{PUL}}$ & mmHg$\cdot$ s/mL & 0.071 & Pulmonary arterial resistance\\
$C_{\mathrm{AR}}^{\mathrm{PUL}}$ & mL/mmHg & 6.0 & Pulmonary arterial compliance\\
$L_{\mathrm{AR}}^{\mathrm{PUL}}$ & mmHg$\cdot$ s$^\mathrm{2}$/mL & 0.000021 & Pulmonary arterial inertia\\
$R_{\mathrm{C}}^{\mathrm{PUL}}$ & mmHg$\cdot$ s/mL & 0.018 & Oxygenated pulmonary capillary resistance\\
$C_{\mathrm{C}}^{\mathrm{PUL}}$ & mL/mmHg & 5.78 & Oxigenated pulmonary capillary compliance\\
$R_{\mathrm{SH}}$ & mmHg$\cdot$ s/mL & 0.35 & Non-oxygenated pulmonary capillary resistance\\
$C_{\mathrm{SH}}$ & mL/mmHg & 0.049 & Non-oxygenated pulmonary capillary compliance\\
$R_{\mathrm{VEN}}^{\mathrm{PUL}}$ & mmHg$\cdot$ s/mL & 0.038 & Pulmonary venous resistance\\
$C_{\mathrm{VEN}}^{\mathrm{PUL}}$ & mL/mmHg & 13.18 & Pulmonary venous compliance\\
$L_{\mathrm{VEN}}^{\mathrm{PUL}}$ & mmHg$\cdot$ s$^\mathrm{2}$/mL & 0.000021 & Pulmonary venous inertia\\
$HR$ & s & $80$ & Heart rate\\
$tC_{\mathrm{LA}}$ & s & $0.75T$$_{\text{HB}}$ & Time of left atrial contraction\\
$TC_{\mathrm{LA}}$ & s & $0.1$T$_{\text{HB}}$ & Duration of left atrial contraction\\
$tR_{\mathrm{LA}}$ & s & $tC_{\mathrm{LA}}+TC_{\mathrm{LA}}$ & Time of left atrial relaxation\\
$TR_{\mathrm{LA}}$ & s & $0.8$T$_{\text{HB}}$ & Duration of left atrial relaxation\\
$tC_{\mathrm{LV}}$ & s & $0.0$ & Time of left ventricular contraction\\
$TC_{\mathrm{LV}}$ & s & $0.265$T$_{\text{HB}}$ & Duration of left ventricular contraction\\
$tR_{\mathrm{LV}}$ & s & $tC_{\mathrm{LV}}+TC_{\mathrm{LV}}$ & Time of left ventricular relaxation\\
$TR_{\mathrm{LV}}$ & s & $0.4T$T$_{\text{HB}}$ & Duration of left ventricular relaxation\\
$tC_{\mathrm{RA}}$ & s & $0.8$T$_{\text{HB}}$ & Time of right atrial contraction\\
$TC_{\mathrm{RA}}$ & s & $0.1$T$_{\text{HB}}$ & Duration of right atrial contraction\\
$tR_{\mathrm{RA}}$ & s & $tC_{\mathrm{RA}}+TC_{\mathrm{RA}}$ & Time of left atrial relaxation\\
$TR_{\mathrm{RA}}$ & s & $0.7$T$_{\text{HB}}$ & Duration of left atrial relaxation\\
$tC_{\mathrm{RV}}$ & s & $0.0$ & Time of right ventricular contraction\\
$TC_{\mathrm{RV}}$ & s & $0.3$T$_{\text{HB}}$ & Duration of left ventricular contraction\\
$tR_{\mathrm{RV}}$ & s & $tC_{\mathrm{RV}}+TC_{\mathrm{RV}}$ & Time of right ventricular relaxation\\
$TR_{\mathrm{RV}}$ & s & $0.4$T$_{\text{HB}}$ & Duration of left ventricular relaxation\\
\hline
\end{tabular}
\caption{List of parameters and their reference values for an ideal healthy individual}
\label{table:params}
\end{center}
\end{table}

\begin{table}[H]
\scriptsize
\makebox[\textwidth][c]{
	\begin{tabular}{|c|c|c|c|c|c|}
	\hline
	&Model output & Unit & Range & Model value & Description\\
	\hline
	\multirow{8}{*}{\begin{tabular}{c}Model outputs\\ used for calibration\end{tabular}} & $LA_{\mathrm{I-Vmax}}$ & mL/m$^2$ & [16,34]\cite{lang2015recommendations} & 22.2 & Indexed maximal left atrial volume\\
	& $LV_{\mathrm{I-EDV}}$ & mL/m$^2$ & [50,90]\cite{hurst1990heart} & 59.7 & Indexed left ventricular end diastolic volume\\
	& $LV_{\mathrm{ESV}}$ &  mL & [18,52]\cite{lang2015recommendations} & 42.7 & Left ventricular end systolic volume\\
	& $LV_{\mathrm{EF}}$ & \% & [53,73]\cite{lang2015recommendations} & 60.0 & Left ventricular ejection fraction\\
	& max $\nabla P_{\mathrm{rAV}}$ & mmHg& - & 17.4 & Maximal right atrioventricular pressure gradient\\
	& $SAP_{\mathrm{max}}$ &  mmHg & [-,140]\cite{lang2015recommendations} & 109.6 & Systolic systemic arterial pressure\\
	& $SAP_{\mathrm{min}}$ &  mmHg & [-,80]\cite{lang2015recommendations} & 71.3 & Diastolic systemic arterial pressure\\
	& $PAP_{\mathrm{max}}$ & mmHg & [15,28]\cite{hurst1990heart} & 23.6 & Systolic pulmonary arterial pressure\\
	\hline	
	\multirow{23}{*}{\begin{tabular}{c}Additional model \\ outputs\end{tabular}} 
	& $LA_\mathrm{Pmax}$ & mmHg & [6,20]\cite{hurst1990heart} & 10.3 & Maximal left atrial pressure\\
	& $LA_\mathrm{Pmin}$ & mmHg & [-2,9]\cite{hurst1990heart} & 5.7 & Minimal left atrial pressure\\
	& $LA_\mathrm{Pmean}$ & mmHg & [4,12]\cite{hurst1990heart} & 8.8 & Mean left atrial pressure\\
	& $LV_{\mathrm{SV}}$ & mL & [30,80]\cite{lang2015recommendations} & 64.1 & Left ventricular stroke volume\\
	& $CI$ & L/min/m$^2$ & [2.8,4.2]\cite{hurst1990heart} & 2.9 & Cardiac index\\
	& $LV_{\mathrm{Pmax}}$ & mmHg & [90,140]\cite{hurst1990heart} & 110.5 & Maximal left ventricular pressure\\
	& $LV_{\mathrm{Pmin}}$ & mmHg & [4,12]\cite{hurst1990heart} & 4.0 & Minimal left ventricular pressure\\
	& $RA_{\mathrm{I-Vmax}}$ & mL/m$^2$ & [10,36]\cite{lang2015recommendations} & 29.3 & Indexed maximal right atrial volume\\
	& $RA_\mathrm{Pmax}$ & mmHg & [2,14]\cite{hurst1990heart} & 9.6 & Maximal right atrial pressure\\
	& $RA_\mathrm{Pmin}$ & mmHg & [-2,6]\cite{hurst1990heart} & 4.4 & Minimal right atrial pressure\\
	& $RA_\mathrm{Pmean}$ & mmHg & [-1,8]\cite{hurst1990heart} & 6.9 & Mean right atrial pressure\\
	& $RV_{\mathrm{I-EDV}}$ &  mL/m$^2$ & [44,80]\cite{rudski2010guidelines} & 68.2 & Indexed right ventricular end diastolic volume\\
	& $RV_{\mathrm{I-ESV}}$ &  mL/m$^2$ & [19,46]\cite{rudski2010guidelines} & 32.6 & Indexed right ventricular end systolic volume\\
	& $RV_{\mathrm{EF}}$ & \% & [44,71]\cite{rudski2010guidelines} & 52.2 & Right ventricular ejection fraction\\
	& $RV_{\mathrm{Pmax}}$ & mmHg & [15,28]\cite{hurst1990heart} & 25.2 & Maximal right ventricular pressure\\
	& $RV_{\mathrm{Pmin}}$ & mmHg & [0,8]\cite{hurst1990heart} & 3.4 & Minimal right ventricular pressure\\
	& $PAP_{\mathrm{min}}$ & mmHg & [5,16]\cite{hurst1990heart} & 15.9 & Diastolic pulmonary arterial pressure\\
	& $PAP_{\mathrm{mean}}$ & mmHg & [10,22]\cite{hurst1990heart} & 19.5 & Mean pulmonary arterial pressure\\
	& $PWP_{\mathrm{min}}$ & mmHg & [1,12]\cite{hurst1990heart} & 11.5 & Minimal pulmonary wedge pressure\\
	& $PWP_{\mathrm{mean}}$ & mmHg & [6,15]\cite{hurst1990heart} & 12.0 &Mean pulmonary wedge pressure \\
	& SVR & mmHg$\cdot$ min/L & [11.3,17.5]\cite{hurst1990heart} & 16.2 & Systemic vascular resistance\\
	& PVR & mmHg$\cdot$ min/L & [1.9,3.1]\cite{hurst1990heart} & 2.09 & Pulmonary vascular resistance\\
	& Shunt Fraction & \% & [0,5]\cite{velthuis2015pulmonary} & 4.73 & Shunt fraction\\
	\hline
	\end{tabular}}

\caption{List of model outputs, the units of measure, the echocardiographic ranges for a healthy individual and the values returned by the numerical model with the reference setting of parameters.}
\label{table:qoi}
\end{table}

\begin{table}[H]
\scriptsize
\begin{center}
	\begin{tabular}{|c|c|c|c|c|c|c|c|}
	\hline
	$LA_\mathrm{Vmax}$ & $LV_{\mathrm{EDV}}$ & $LV_{\mathrm{ESV}}$ & $LV_{\mathrm{EF}}$ & max $\nabla P_{\mathrm{rAV}}$ & $SAP_{\mathrm{max}}$ 	&  $SAP_{\mathrm{min}}$ & $PAP_{\mathrm{max}}$ \\
	\hline
	$5\%$ & $5\%$ & $5\%$ & $4\%$ & $4\%$ & $4\%$ & $5\%$ & $5\%$\\
	\hline
	\end{tabular}
\end{center}
\caption{Measurement errors on the model outputs used as standard deviations for generating the noisy data.}
\label{table:measErr}
\end{table}

\begin{table}[H]
\scriptsize
\begin{center}
	\begin{tabular}{|c|c|c|c|c|c|c|c|c|c|}
	\hline
	Patient & HR & $LA_\mathrm{Vmax}$ & $LV_{\mathrm{EDV}}$ & $LV_{\mathrm{ESV}}$ & $LV_{\mathrm{EF}}$ & max $\nabla P_{\mathrm{rAV}}$ & $SAP_{\mathrm{max}}$ 	&  $SAP_{\mathrm{min}}$ & $PAP_{\mathrm{max}}$ \\
	\hline
	Monzino & $70$ & - & $233$ & $130$ & $42$ & $25$ & $140$ & $55$ & - \\
	Sacco & $60$ & $50$ & $110$ & $33$ & $70$ & $20$ & $135$ & $70$ & $25$ \\
	\hline
	\end{tabular}
\end{center}
\caption{Clinical data of two patients, provided by Centro Cardiologico Monzino and L. Sacco Hospital in Milan.}
\label{table:real_test}
\end{table}

\begin{table}[H]
\scriptsize
\begin{center}
	\begin{tabular}{|c|c|c|c|}
	\hline
	Patient & $RMSE$ CMC & $RMSE$ L-BFGS-B & $RMSE$ CMC-L-BFGS-B\\
	\hline
	Monzino & $5.4\cdot 10^{-2}$ & $4.0\cdot 10^{-2}$ & $3.8\cdot 10^{-2}$\\
	Sacco & $4.0\cdot 10^{-2}$ & $1.2\cdot 10^{-2}$ & $1.2\cdot 10^{-2}$\\
	\hline
	\end{tabular}
\end{center}
\caption{Final value of the $RMSE$ for the three patient-specific calibration methods.}
\label{table:losses}
\end{table}

\begin{table}[H]
\scriptsize
\begin{center}
	\begin{tabular}{|c|c|c|c|c|c|c|c|c|}
	\hline
	Patient & Calibration method & $EB_\mathrm{LA}$ & $EA_\mathrm{LV}$ & $EB_\mathrm{LV}$ & $EA_\mathrm{RV}$ & $R_\mathrm{AR}^\mathrm{SYS}$ & $C_\mathrm{AR}^\mathrm{SYS}$ & $R_\mathrm{VEN}^\mathrm{SYS}$\\
	\hline
	\multirow{3}{*}{\begin{tabular}{c}Monzino\end{tabular}}& CMC & - & $1.04$ & - & $0.66$ & $0.62$ & $0.65$ & $0.27$ \\
	& L-BFGS-B  & - & $1.02$ & - & $0.66$ & $0.55$ & $0.75$ & $0.24$\\
	& CMC-L-BFGS-B & - & $1.02$ & - & $0.60$ & $0.55$ & $0.75$ & $0.23$ \\
	Relative standard deviation & & - & $1\%$ & - & $5\%$ & $7\%$ & $8\%$ & $8\%$ \\
	\hline
	\multirow{3}{*}{\begin{tabular}{c}Sacco\end{tabular}} & CMC & $0.21$ & $4.18$ & $0.067$ & $0.62$ & $0.63$ & $1.08$ & $0.301$  \\
	& L-BFGS-B & $0.18$ & $4.38$ & $0.061$ & $0.72$ & $0.60$ & $0.89$ & $0.299$\\
	& CMC-L-BFGS-B & $0.19$ & $4.35$ & $0.061$ & $0.72$ & $0.60$ & $0.90$ & $0.297$ \\
	Relative standard deviation & & $7\%$ & $2\%$ & $6\%$ & $8\%$ & $3\%$ & $11\%$ & $1\%$ \\
	\hline
	\end{tabular}
\end{center}
\caption{Estimated parameters by the patient-specific calibrations and their relative standard deviations. For the Monzino patient, we do not calibrate $EB_\mathrm{LA}$ and $EB_\mathrm{LV}$ according to the global sensitivity analysis.}
\label{table:est_params}
\end{table}

\newpage

\section*{Figures}

\end{document}